\documentclass[12pt]{article}
\pdfoutput=1
\usepackage[utf8]{inputenc}
\usepackage{fullpage}
\usepackage[small,bf]{caption}
\usepackage{fancyvrb}
\usepackage{float}
\usepackage{geometry}
\usepackage{verbatim}
\usepackage{longtable}
\usepackage{graphicx}
\usepackage{enumerate, enumitem}
\usepackage{booktabs}
\usepackage[
colorlinks=true,
citecolor=blue,
linkcolor=magenta,
urlcolor=cyan]{hyperref}
\usepackage{url}
\usepackage{xcolor}
\usepackage{tikz}
\usepackage{xparse}
\usepackage{amsmath,amsfonts}
\usepackage{subcaption}
\usepackage[american]{circuitikz}

\newcommand{\eg}{{\it e.g.}}
\newcommand{\ie}{{\it i.e.}}

\newcommand{\BEQ}{\begin{equation}}
\newcommand{\EEQ}{\end{equation}}
\newcommand{\BEAS}{\begin{eqnarray*}}
\newcommand{\EEAS}{\end{eqnarray*}}
\newcommand{\ones}{\mathbf 1}
\newcommand{\reals}{{\mbox{\bf R}}}

\newcommand{\complex}{{\mbox{\bf C}}}
\newcommand{\symm}{{\mbox{\bf S}}}  % symmetric matrices
\newcommand{\sign}{{\mbox{\bf sign}}}

\newcommand{\Tr}{\mathop{\bf Tr}}
\newcommand{\diag}{\mathop{\bf diag}}

\newcommand{\Expect}{\mathop{\bf E{}}}

\newcommand{\dom}{\mathop{\bf dom}}

\usepackage{relsize} % Relsize package needed

\usepackage{booktabs}
\usepackage{array}

\usepackage{listings}

\definecolor{codegreen}{rgb}{0,0.6,0}
\definecolor{codegray}{rgb}{0.5,0.5,0.5}
\definecolor{codepurple}{rgb}{0.58,0,0.82}
\definecolor{backcolour}{rgb}{0.95,0.95,0.92}

\lstdefinestyle{mystyle}{
    backgroundcolor=\color{backcolour},   
    commentstyle=\color{codegreen},
    keywordstyle=\color{magenta},
    numberstyle=\tiny\color{codegray},
    stringstyle=\color{codepurple},
    basicstyle=\ttfamily\small,
    breakatwhitespace=false,         
    breaklines=true,                 
    captionpos=b,                    
    keepspaces=true,                 
    numbers=left,                    
    numbersep=5pt,                  
    showspaces=false,                
    showstringspaces=false,
    showtabs=false,                  
    tabsize=2
}

\definecolor{seagreen}{rgb}{0.18, 0.55, 0.34}
\definecolor{mediumviolet-red}{rgb}{0.78, 0.08, 0.52}
\definecolor{khaki}{rgb}{0.94, 0.9, 0.55}

\lstdefinelanguage{mypython}
{
	keywords=[1]{from, import, assert, not, print},
	keywordstyle=[1]{\color{mediumviolet-red}},
	keywords=[2]{surecr, torch, cp, lo, pl},
	keywordstyle=[2]{\color{seagreen}},
	numbers=none,
	upquote=true,
	showstringspaces=false,
	basicstyle=\ttfamily,
	columns=fullflexible,
	keepspaces=true,
	emph={True,False,as,def,return,float,class,match,switch,len},
	emphstyle={\color{seagreen}},
	frame=trBL,
	belowskip=1em,
	aboveskip=1em,
	captionpos=b
}

\usetikzlibrary{math}

\usepackage{textcomp}

\begin{document}

\title{Disciplined Nonlinear Programming}
\author{Daniel Cederberg \and William Zhang \and Parth Nobel \and Stephen Boyd}

\date{\today}

\maketitle

\begin{abstract}
We introduce \emph{disciplined nonlinear programming} (DNLP), a syntax for
specifying nonlinear programming problems. DNLP is inspired by disciplined
convex programming (DCP) and allows smooth functions to be freely mixed with
nonsmooth convex and concave functions, with rules governing how the nonsmooth
functions can be used. Problems expressed in DNLP form can be automatically
canonicalized to a standard nonlinear programming (NLP) form and passed to a
suitable NLP solver.  As in DCP, the canonicalization relaxes nonsmooth convex
and concave functions in a lossless way, allowing them to be handled by NLP
solvers that require smooth functions. In addition to extending NLP to include
useful nondifferentiable convex and concave functions, transforming the original
problem to an equivalent NLP form offers several advantages, including simpler
problem initialization. We describe the language and our open-source
implementation of DNLP as an extension of CVXPY, a parser for DCP.
\end{abstract}

\clearpage
\tableofcontents
\clearpage

\section{Introduction}
\subsection{Nonlinear programming}
Nonlinear programming (NLP) has a long and well-established history
\cite{giorgi2013traces}, with successful applications spanning decades in fields
such as chemical engineering \cite{biegler2010nonlinear}, topology optimization
\cite{bendsoe_topology_2004}, optimal control \cite{betts2010practical},
aerospace design \cite{keane2005computational}, and design optimization
\cite{PapalambrosWilde2017}, among others. This breadth of applications
highlights the remarkable generality of NLP as a unifying framework for modeling
and solving problems.

However, the generality of NLP comes at a cost. With the exception of global
optimization methods \cite{Horst1996}, which are often computationally
prohibitive, there are no universal guarantees of achieving global optimality,
and in many cases solving NLPs remains as much an art as a science. While the
usual concern is the lack of global optimality guarantees, other pathologies can
occur, including failure to converge to a feasible point even when one exists.
NLP solvers will do their best to find a solution, but success depends on how
the problem is formulated, the choice of algorithm, its hyperparameters, and the
initialization. Nonetheless, NLP remains a powerful and widely used tool, as
evidenced by the popularity of general-purpose NLP solvers such as Ipopt
\cite{wachter2006}.

\subsection{Modeling languages}
To interface with NLP solvers, several modeling languages have been developed.
Classic examples include the commercial systems AMPL \cite{fourer1990modeling},
AIMMS \cite{bisschop2006aimms}, and GAMS \cite{brook1988gams}, which are based
on their own domain-specific programming languages. More recent open-source
frameworks are instead embedded in general-purpose languages, such as YALMIP
\cite{lofberg2004yalmip} in MATLAB, JuMP \cite{dunning2017jump} in Julia, Pyomo
\cite{hart2011pyomo, bynum2021pyomo} in Python, and CasADi
\cite{andersson2019casadi} in C++. (We note that AMPL, GAMS, and CasADi all have
MATLAB and Python interfaces as well.) These modeling languages facilitate the
specification of NLPs but largely treat user-specified problem formulations as
black boxes. As a result, a poorly structured formulation may be passed to the
solver, making it difficult for the solver to find a solution. (We give two such
examples in \S \ref{sec:advantages-of-canonicalization}.)

In this paper, we take the viewpoint that an NLP modeling language should (to
the extent possible) exploit the structure of the user-specified problem
formulation and reformulate it to increase the likelihood that the solver
succeeds. To this end, we introduce a grammar for specifying NLPs, which we call
\emph{disciplined nonlinear programming} (DNLP). To handle nonsmooth convex and
concave functions, DNLP adopts the same core idea as \emph{disciplined convex
programming} (DCP) \cite{grant2004disciplined, grant2006disciplined}, a grammar
for specifying convex optimization problems, and analyzes monotonicity to relax
nonsmooth functions into equivalent smooth formulations
\cite{graphimplementations2008}. The popular convex optimization modeling
language CVXPY \cite{diamond2016cvxpy, agrawal2018rewriting} is based on DCP,
and we have implemented a rewriting system based on DNLP as an extension to
CVXPY. This extension allows users to seamlessly specify NLPs as long as they
conform to a minimal set of rules, and the problem is then (hopefully) solved by
an NLP solver.

It is important to note, however, that the discipline imposed by DNLP does not,
in itself, guarantee that a solver will succeed and be able to compute a
solution. But we believe that following the DNLP ruleset increases the
likelihood of successful convergence. This should be contrasted with convex
optimization and DCP, where the benefits of imposing such discipline are much
stronger: any formulation conforming to DCP is automatically certified as convex
and can be solved reliably and efficiently to global optimality (up to some
practical problem size limits and solver tolerances).

\subsection{Outline}
The remainder of this paper begins with a brief overview of NLP in \S
\ref{sec:nlp}. In \S \ref{sec:dnlp}, we introduce DNLP and its (minimal)
ruleset, while \S \ref{sec:canonicalization} describes the canonicalization
process and explains how DNLP allows nonsmooth problems to be relaxed (without
loss) into equivalent smooth formulations. Finally, in \S
\ref{sec:planning-geometry}-\S\ref{sec:statistics}, we present several numerical
examples from many different fields.

\section{Nonlinear programming}
\label{sec:nlp}
A \emph{nonlinear program} is an optimization problem of the form
\begin{equation} \label{e:standard-form}
\begin{array}{ll}
\mbox{minimize} & f(x) \\
\mbox{subject to} & c(x) = 0  \\
& \ell \leq x \leq u,  \\
\end{array}
\end{equation}
or one that can be readily converted into this form. Here, $x \in \reals^n$ is
the optimization variable, $\ell \in \reals^n$ and $u \in \reals^n$ are given
variable bounds, and $f: \reals^n \to \reals$ and $c: \reals^n \to \reals^m$ are
differentiable functions that are allowed to be nonconvex. An inequality
constraint of the form $d_i(x) \leq 0$ can be expressed in this form by
introducing a slack variable $s_i \geq 0$ together with the constraint $d_i(x) +
s_i = 0$. An unbounded variable $x_i$ can be specified by setting $\ell_i =
-\infty$ and $u_i = \infty$.

In this section we provide a survey of NLP, including common variations on the
standard form given above, algorithms and solvers, and theoretical properties.
For more background we refer the reader to the many excellent textbooks on the
subject \cite{nocedal2006numerical, Fletcher2000, gill2019practical,
bertsekas2016}.

\subsection{Standard forms and oracles}
Many NLP solvers have been developed over the years (we name a few of these in
\S \ref{sec:solvers}), each with its own interface and its own standard form.
While different solvers have their own standard forms, they are all closely
related to \eqref{e:standard-form}, or they convert problems into this form
internally. For example, Ipopt \cite{wachter2006} requires constraints to be
given as two-sided inequalities of the form $\ell \leq g(x) \leq u$, and
internally transforms the constraint into the form \eqref{e:standard-form} by
introducing a new variable $s$ together with the equality constraint $g(x) - s =
0$ and the bounds $\ell \leq s \leq u$. Other solvers, such as Knitro
\cite{byrd2006} or SNOPT \cite{gill2005}, allow users to specify linear
constraints separately for further efficiency. Manually reformulating an
optimization problem to match a solver’s standard form is tedious and prone to
errors. Modeling languages automate this process, allowing users to switch
seamlessly between solvers with different standard forms.

In addition to transforming user-specified problems into the standard form
expected by solvers, NLP modeling languages are responsible for providing
oracles that evaluate the objective and constraint functions and their
derivatives. Most, if not all, modeling languages construct these oracles
using \emph{automatic differentiation} \cite{griewank2008}.

\subsection{Algorithms and solvers}
\label{sec:solvers}
Algorithms for solving NLPs have been studied since at least the 1940s (see,
\eg, \cite{giorgi2013traces}), but only in the past few decades, with advances
in software, have these methods become accessible to a broader audience. The two
most common types of algorithms implemented in modern NLP solvers are
\emph{interior-point methods} (IPMs) and \emph{sequential quadratic programming}
(SQP).

Interior-point methods reduce \eqref{e:standard-form} to a sequence of
equality-constrained problems by incorporating the inequality constraints into
the objective using \emph{barrier functions}. A large body of theory on barrier
functions for solving NLPs was developed during the 1960s \cite{fiacco1968}, but
researchers lost interest in the most basic IPM—the \emph{primal barrier
method}—due to concerns about ill-conditioning \cite{murray1971} that later
proved unfounded \cite{wright1998, Forsgren2002}. Much later, more sophisticated
IPMs for NLP were developed, and today many of the most popular solvers
implement IPMs, including the open-source solvers Ipopt \cite{wachter2006} and
Uno \cite{Vanaret2025}, as well as the commercial solvers LOQO
\cite{vanderbei1999}, Knitro \cite{byrd2006}, and Gurobi \cite{gurobiVersion13}.

Sequential quadratic programming methods reduce \eqref{e:standard-form} to a
sequence of quadratic programs. The constraints of each quadratic subproblem are
linearizations of the constraints in the original problem, and the objective is
a quadratic approximation of the Lagrangian function. SQP methods were first
proposed in the 1960s \cite{wilson1963}, and modern solvers implementing SQP
include the commercial packages SNOPT \cite{gill2005}, Knitro-Active
\cite{byrd2006}, and WORHP \cite{buskens2012esa}, as well as the open-source
solvers GRANSO \cite{Curtis2017, liang2022ncvx} and Uno \cite{Vanaret2025}.

While IPMs and SQP methods are the most commonly implemented algorithms, several
solvers also implement \emph{augmented Lagrangian methods}, including the
open-source solver Algencan \cite{andreani2008} and the commercial solvers MINOS
\cite{Murtagh1983}, Lancelot \cite{Conn1992}, and the recent Knitro-Augmented
\cite{knitroAugmented}. These methods reduce \eqref{e:standard-form} to a
sequence of subproblems in which the objective is the Lagrangian augmented with
a penalty term for constraint violations, where some or all of the constraints
are incorporated into the penalty and the remaining constraints are enforced
explicitly.

Given the many NLP solvers available, a natural question is which solver to use
for a given problem. While there is no definitive answer, the conventional
wisdom is that IPMs are faster and more reliable when solving a problem from
scratch, \ie, without a good initial point \cite{gill2015}. However, IPMs such
as Ipopt may struggle with problems that violate standard regularity conditions
(see, \eg, \cite[\S 11]{biegler2010nonlinear} or \cite{Thierry2020}), in which
case augmented Lagrangian methods can be more robust \cite{izmailov2012}. For
example, Knitro states on their website that the primary advantage of their
augmented Lagrangian method over IPMs is that it is ``designed to better handle
difficult problems with degenerate constraints where the linear independence
constraint qualification (LICQ) is not satisfied''. Nevertheless, we recommend
trying Ipopt first, because it is open-source, widely adopted (as evidenced by
its citation count), and performs well across many applications. In our
experience, it works very well.

\subsection{Theoretical properties}
\begin{figure}
  \centering
  \begin{subfigure}[b]{\textwidth}
    \centering
    \includegraphics[width=0.55\textwidth]{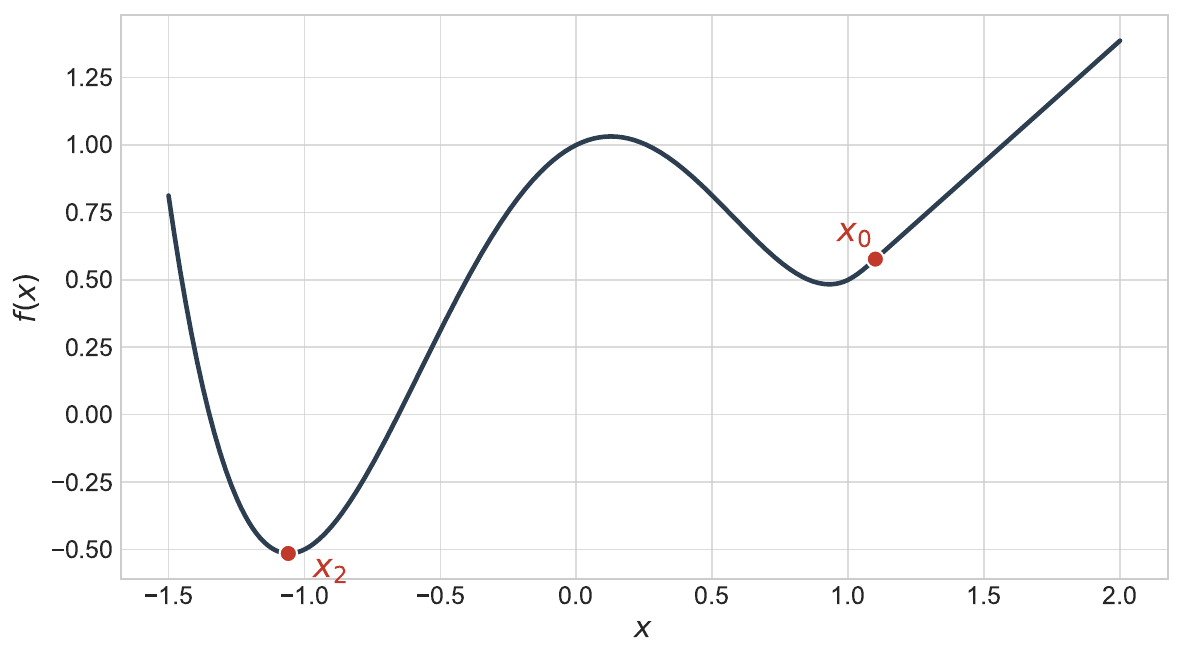}
    \caption{}
    \label{fig:local_myth}
  \end{subfigure}
  \begin{subfigure}[b]{\textwidth}
    \centering
    \includegraphics[width=0.55\textwidth]{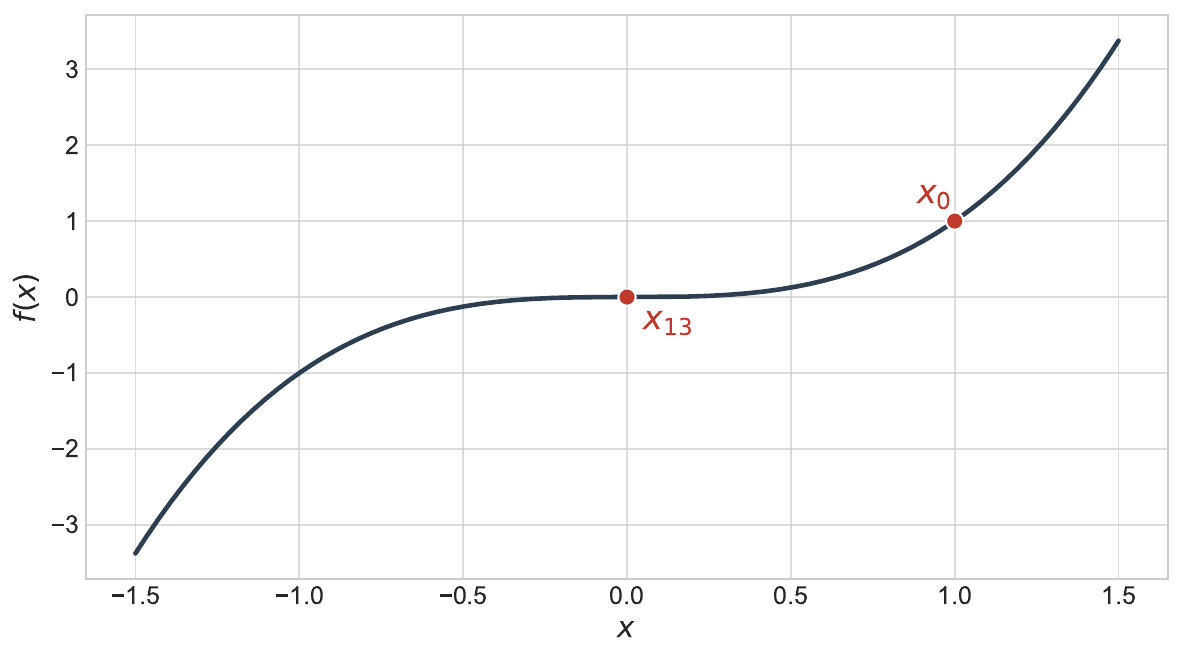}
    \caption{}
    \label{fig:saddle_point}
  \end{subfigure}
  \begin{subfigure}[b]{\textwidth}
    \centering
    \includegraphics[width=0.55\textwidth]{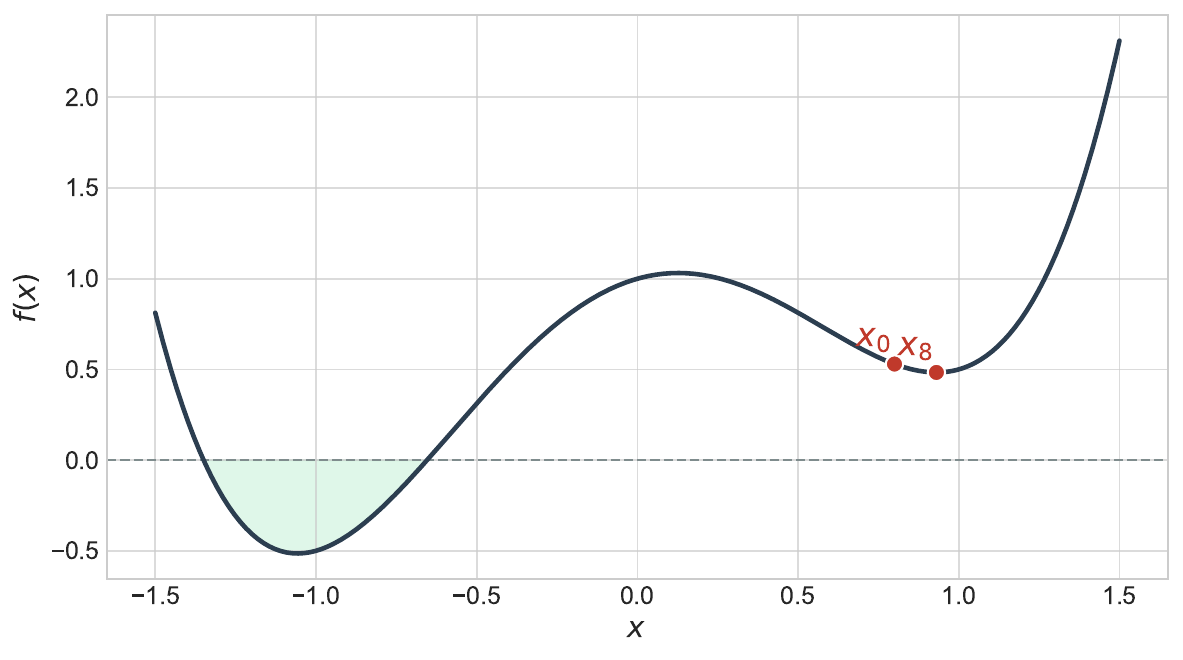}
    \caption{}
    \label{fig:infeas}
  \end{subfigure}
  \caption{The red dots denote the initial point and the iterate that Ipopt converged to.
  \textbf{(a)} Ipopt converges (in two iterations) to a local minimizer that is not the
  nearest one to the initial point. 
  \textbf{(b)} Ipopt converges (in thirteen iterations) to a saddle point. 
  \textbf{(c)} Ipopt converges (in eight iterations) to an infeasible point. }
  \label{fig:convergence}
\end{figure}

Because NLP covers a vast range of problems, including many known to be NP-hard,
it is unrealistic to expect NLP solvers to guarantee convergence to a
\emph{global} minimizer (\ie, a feasible point achieving the smallest possible
objective value among all feasible points). While this misconception is widely
recognized, others are more subtle. For example, each of the following
statements is false, even though each is weaker than the one preceding it:
\begin{itemize}
\item NLP solvers always converge to the local minimizer nearest the initial point.
\item NLP solvers always converge to a local minimizer.
\item NLP solvers always converge to a point that approximately satisfies a set of
necessary but not sufficient optimality equations known as the
\emph{Karush-Kuhn-Tucker} (KKT) conditions \cite[\S 4]{bertsekas2016}.
\end{itemize}
The first two statements are false even for unconstrained problems, and the
third is false because a solver may fail to converge to a feasible point even
when one exists. We present counterexamples to all three statements in figure
\ref{fig:convergence}. (A similar example to the third case is given in
\cite{Wachter2026}.) In the top panel, we minimize a twice-continuously
differentiable function that is very flat at the initial point, so Ipopt steps
over the nearest local minimizer and converges to a different one. In the middle
panel, we minimize $f(x) = x^3$ starting at $x_0 = 1$ and Ipopt converges to the
saddle point at $x=0$. In the bottom panel, we minimize $f(x) = x^4 - 2x^2 +
0.5x + 1$ subject to $f(x) \leq 0$. The feasible region is shaded in green, and
Ipopt converges to an infeasible point.

Despite these counterexamples, NLP solvers in practice often successfully compute
KKT points, \ie, points satisfying the KKT conditions within some tolerances.
When an NLP solver claims that it has solved a nonlinear program, it typically
means that it has found a KKT point. However, not all KKT points are local
minimizers, so solvers incorporate various heuristics and techniques to steer
iterates away from such undesirable KKT points.

While not every KKT point is a local minimizer, the converse is true under
so-called \emph{constraint qualifications} \cite[\S 4.3.4]{bertsekas2016}. The
derivation of the KKT conditions is based on the idea of linearizing the
constraints around a local minimizer, and constraint qualifications are
conditions that ensure that this linearization is a good approximation of the
true constraints. A common constraint qualification assumed by NLP solvers is
the \emph{linear independence constraint qualification} (LICQ), which for
problem \eqref{e:standard-form} requires that the set consisting of the
gradients of the active bound constraints and the gradients of the equality
constraints is linearly independent. If LICQ holds at a local minimizer
$x^\star$, then $x^\star$ is guaranteed to also be a KKT point, and the
so-called \emph{Lagrange multipliers} (which are auxilliary variables in the KKT
conditions) are guaranteed to be unique. LICQ also seems to play a role in
practice; if it does not hold at a local minimizer, some solvers are less
robust and may fail to converge.

\section{Disciplined nonlinear programming}
\label{sec:dnlp}
Nondifferentiable functions often arise in applications and pose significant
challenges for most NLP solvers. A na\"{i}ve approach is to simply ignore these
nondifferentiabilities or assume they will not occur in practice, but this
often leads to poor performance and solver failure. The difficulty is that
the points of nondifferentiability are often precisely the points of interest.
For example, in problems with $\ell_1$ regularization, the goal
is typically to find a solution in which the argument is sparse,
which is a point where the $\ell_1$-norm is nondifferentiable.

To support nondifferentiable functions in nonlinear programs without
compromising solver reliability, we introduce the notion of \emph{disciplined
nonlinear programming} (DNLP). It consists of two key components:
\begin{itemize}
\item An \emph{atom library}---a collection of
functions that can be used to describe a problem.
These functions have known attributes including smoothness, sign, monotonicity,
and curvature.
\item The \emph{DNLP ruleset}---a set of rules specifying how atoms may be
combined to form more complicated expressions, and how these expressions may appear
in objectives and constraints.
\end{itemize}
This framework guarantees that any problem with nonsmooth functions complying
with the DNLP ruleset admits an equivalent formulation in which all the involved
functions are smooth. Furthermore, the canonicalization process preserves
standard NLP regularity conditions, such as LICQ (\ie, if the original problem
satisfies LICQ at a local minimizer, then the canonicalized problem also
satisfies LICQ at the corresponding local minimizer.) For problems that only
involve smooth functions, DNLP imposes no additional restrictions.

DNLP is heavily inspired by DCP, so readers familiar with DCP
will find many similarities. Roughly speaking, DNLP mirrors the structure of
DCP, with smooth functions playing the role of affine functions, and
generalizations of convex and concave functions that can be mixed with smooth
functions.

\subsection{Atoms}
The rules of DNLP depend on the smoothness and curvature properties of atoms. We
classify atoms into three categories: \emph{smooth}, \emph{nonsmooth-convex}
(NS-convex), and \emph{nonsmooth-concave} (NS-concave). We list some atoms and
their classifications in table \ref{tab:smoothness-classification}.

\paragraph{Smooth atoms.}
An atom is \emph{smooth} if it is twice continuously differentiable in the
\emph{interior} of its domain. For example, the atoms $\phi_{\text{log}}$ and
$\phi_{\text{sqrt}}$ defined by $\phi_{\text{log}}(x) = \log x$ with $\dom
\phi_{\text{log}} = \reals_{++}$, and $\phi_{\text{sqrt}}(x) = \sqrt{x}$ with $\dom
\phi_{\text{sqrt}} = \reals_{+},$ are both smooth, and so is any
affine or trigonometric atom. In contrast, the atom $\phi_{\text{abs}}$ defined
by $\phi_{\text{abs}}(x) = |x|$ with $\dom \phi_{\text{abs}} = \reals$ is not
smooth. 

\paragraph{Nonsmooth-convex atoms.}
An atom is \emph{NS-convex} if it is convex and not twice continuously
differentiable in the interior of its domain. Two examples are 
$\phi_{\text{max}}$ and $\phi_{\text{norm2}}$ defined by
$\phi_{\text{max}}(x, y) = \max(x, y)$ and $\phi_{\text{norm2}}(x) = \| x \|_2$.
(The latter is not differentiable at $x=0$.)

\paragraph{Nonsmooth-concave atoms.}
An atom is \emph{NS-concave} if it is concave and not twice continuously differentiable
in the interior of its domain. Two examples are $\phi_{\text{min}}$ and
$\phi_{\text{sum\_smallest}}$ defined by $\phi_{\text{min}}(x, y) = \min(x, y)$
and $\phi_{\text{sum\_smallest}}(x; k) = \sum_{i=n-k+1}^n x_{[i]}$, where $x_{[i]}$
is the $i$th largest element of $x \in \reals^n$, and $k \in \{1, \ldots, n\}$
is a fixed parameter.

\begin{table}
\centering
\scriptsize 
\caption{Some atoms and their classifications. If the domain of an atom is not
specified, it means that the atom has full domain.}
\label{tab:smoothness-classification}
\renewcommand{\arraystretch}{1.2} % increase row height
\setlength{\tabcolsep}{2pt} % adjust column spacing
\begin{tabular}{>{\raggedright\arraybackslash}p{2.8cm}
>{\raggedright\arraybackslash}p{7.2cm}
>{\raggedright\arraybackslash}p{2.5cm}
%>{\raggedright\arraybackslash}p{2cm}
}
\toprule
\textbf{Atom} & \textbf{Definition} & \textbf{Domain} \\
\midrule
\multicolumn{3}{c}{\scriptsize{\textbf{Smooth, nonconvex and nonconcave}}} \\
\midrule
\verb|multiply| & $\phi(x, y) = x y$ &  \\ 
\verb|matmul| & $\phi(X, Y) = XY $  \\ 
%\verb|prod & $\phi(x) = \prod_{i=1}^n x_i $ & $x \in \reals^{n}$ \\ 
\verb|quad_form| & $\phi(x) = x^T Q x$ where $Q \in \symm^n$ & \\
\verb|sin| & $\phi(x) = \sin x$ &  \\ 
%\verb|cos| & $\phi(x) = \cos x$ &  \\ 
\verb|tan| & $\phi(x) = \tan x$ & $x \in (-\pi/2, \pi/2)$  \\ 
\verb|sinh| & $\phi(x) = (e^x - e^{-x})/2$ & \\ 
\verb|tanh| & $\phi(x) = (e^x - e^{-x})/(e^x + e^{-x})$ & \\ 
\verb|asinh| & $\phi(x) = \log(x + \sqrt{x^2 + 1})$ &\\ 
\verb|atanh| & $\phi(x) = (1/2) \log\left((1+x)/(1-x)\right)$ & $x \in (-1, 1)$ \\ 
\verb|sigmoid| & $\phi(x) = 1 / (1 + e^{-x})$ & \\ 
\verb|normcdf| & $\phi(x) = (1/\sqrt{2\pi}) \int_{-\infty}^x e^{-t^2/2} \; dt$ &  \\
\midrule
\multicolumn{3}{c}{\scriptsize{\textbf{Smooth, convex or concave}}} \\
\midrule
\verb|exp| & $\phi(x) = e^x$ &  \\ 
\verb|log| & $\phi(x) = \log x$ & $x > 0$ \\
%\verb|entr} & $\phi(x) = -x \log x$ & $x > 0$ \\
%\verb|sum\_squares} & $\phi(x) = x^T x$ & $x \in \reals^n$ \\
\verb|log_sum_exp| & $\phi(x) = \log \left(\sum_{i=1}^n e^{x_i}\right)$ &  \\
\verb|power| & $\phi(x) = x^p$ \text{ where $p > 0$ is an integer} &  \\
\verb|power_pos| & $\phi(x) = x^p$ \text{ where $p > 0$} & $x \geq 0$ \\
%\verb|power\_pos}$ & $\phi(x) = x^p$ \text{ where $p < 0$} & $x > 0$ \\
\verb|sqrt| & $\phi(x) = \sqrt{x}$ & $x \geq 0$ \\ 
%\verb|entr$ & $\phi(x) = -x \log x$ & $x > 0$ \\ 
%\verb|geo\_mean}$ & $\phi(x) = \prod_{i=1}^n x_i^{1/n} $ & $x \in \reals^{n}_+$ \\ 
\verb|inv_pos| & $\phi(x) = 1 / x$ & $x > 0$ \\ 
\verb|quad_over_lin| & $\phi(x, y) = x^Tx / y $ & $y > 0$ \\
%\verb|rel\_entr}$ & $\phi(x, y) = x \log(x / y) $ & $x, y > 0$ \\
\midrule
\multicolumn{3}{c}{\scriptsize{\textbf{Nonsmooth, convex}}} \\
\midrule    
\verb|abs| & $\phi(x) = |x|$ &  \\ 
\verb|max| & $\phi(x) = \max \{x_1, x_2, \ldots, x_n\}$ &  \\ 
\verb|norm1| & $\phi(x) = \|x\|_1$ &  \\ 
\verb|norm2| & $\phi(x) = \|x\|_2$ &  \\ 
\verb|norm_inf| & $\phi(x) = \|x\|_\infty$ &  \\ 
\verb|huber| & $\phi(x; M) =
\left\{
\begin{array}{ll}
x^2, & |x| \leq M \\
2M|x| - M^2, & |x| > M,
\end{array}
\right. 
$ where $M \geq 0$ & \\ 
\verb|sum_largest| & $\phi(x; k) = \sum_{i=1}^k x_{[i]}$ where $k \in \{1, \ldots, n\}$ &  \\ 
\midrule
\multicolumn{3}{c}{\scriptsize{\textbf{Nonsmooth, concave}}} \\
\midrule
\verb|min| & $\phi(x) = \min \{x_1, x_2, \ldots, x_n\}$ &  \\ 
\verb|sum_smallest| & $\phi(x; k) = \sum_{i=n-k+1}^k x_{[i]}$ where $k \in \{1, \ldots, n\}$ & \\ 
\bottomrule
\end{tabular}
\end{table}

\paragraph{Additional attributes.}
Functions in the atom library are also characterized by their sign and
monotonicity. Three categories of monotonicity are considered:
\emph{nondecreasing}, \emph{nonincreasing}, and \emph{nonmonotonic}. The usual
mathematical definitions of monotonicity apply. For functions with multiple
arguments, we specify the monotonicity with respect to each argument separately.
Furthermore, we use \emph{sign-dependent} monotonicity, \ie, the monotonicity of
an atom can depend on the signs of its arguments. For example, the atom
\texttt{square} given by $\phi(x) = x^2$ is classified as nondecreasing for $x
\geq 0$.

\subsection{Expressions}
\label{sec:dnlp-expressions}
An \emph{expression} is recursively defined as an atom evaluated at a
\emph{subexpression}. The subexpression can be a variable, a constant, or
another expression itself. Mathematically, an expression is of the form $f(x) =
\phi(g(x))$ where $\phi$ is the atom and $g(x) = (g_1(x), \ldots, g_k(x))$ is
its argument, the subexpression. We classify expressions into the three
categories \emph{smooth}, \emph{linearizable-convex} (L-convex), and
\emph{linearizable-concave} (L-concave).

\paragraph{Smooth expressions.}
An expression $f(x) = \phi(g(x))$ is defined to be \emph{smooth} if both the
atom $\phi$ and the subexpression $g(x)$ are smooth. Constant expressions and
variable expressions are considered as smooth, so any smooth atom $\phi$,
evaluated at variables or constants, is a smooth expression $\phi(x)$.

\paragraph{L-convex expressions.}
An expression $f(x) = \phi(g(x))$ is defined to be \emph{L-convex} if the atom
$\phi$ is smooth or NS-convex, and for each $i=1, \ldots, k$, one of the
following holds: $g_i(x)$ is smooth; or $g_i(x)$ is L-convex and $\phi$ is
nondecreasing in its $i$th argument; or $g_i(x)$ is L-concave and $\phi$ is
nonincreasing in its $i$th argument.

\paragraph{L-concave expressions.}
An expression $f(x) = \phi(g(x))$ is defined to be \emph{L-concave} if the atom
$\phi$ is smooth or NS-concave, and for each $i = 1, \ldots, k$, one of the
following holds: $g_i$ is smooth; or $g_i$ is L-convex and $\phi$ is
nonincreasing in its $i$th argument; or $g_i$ is L-concave and $\phi$ is
nondecreasing in its $i$th argument. 

\paragraph{Simple consequences of the definitions.} We mention that any smooth
expression is also both L-convex and L-concave. Furthermore, the sum of two
L-convex expressions is L-convex, and the sum of two L-concave expressions is
L-concave. (All these statements follow directly from the definitions of
L-convexity and L-concavity.) This logic is analogous to how, in convex
optimization, affine expressions are both convex and concave, and the sum of two
convex (concave) expressions is convex (concave).

\subsection{Objectives and constraints}
\label{sec:dnlp-ruleset}
For an optimization problem to be a \emph{disciplined nonlinear program}, its
objective and constraints must satisfy the following rules. 

\paragraph{Objective.}
A valid objective is either the minimization of an L-convex expression or the
maximization of an L-concave expression. Maximizing an L-convex expression or
minimizing an L-concave expression is not valid (unless the expression is also
smooth).

\paragraph{Constraints.}
A valid constraint is one of the following:
\begin{itemize}
\item An equality constraint between a smooth left-hand side (LHS) and a
smooth right-hand side (RHS).
\item A less-than-or-equal-to inequality with an L-convex LHS and an
L-concave RHS.
\item A greater-than-or-equal-to inequality with an L-concave LHS and an
L-convex RHS.
\end{itemize}
A problem description that conforms to these rules is called \emph{DNLP-compliant}. We
will see that such a problem formulation can be canonicalized to an equivalent (smooth) NLP
without introducing LICQ violations.

\subsection{Examples}
\label{sec:expression-examples}
\paragraph{DNLP expressions.}
We now give a few examples of expressions that conform to the DNLP ruleset, and
others that do not.   
\begin{itemize}
\item 
The function $f(x, y) = x / y$ with $y > 0$ can be expressed as
\[
\verb|multiply(x, inv_pos(y))|.
\]
When expressed this way, $f(x, y)$ is a smooth expression since it is the
composition of the smooth atom \verb|multiply| with two smooth expressions.
(A variable or a smooth atom by itself is considered a smooth expression; see \S
\ref{sec:dnlp-expressions}.)
\item The function $f(x) = c^T x / (x^T A x) $ with $A \in \symm^n_{++}$  
can be expressed as
\[
\verb|multiply(c @ x, inv_pos(quad_form(x, A)))|.
\]
When expressed this way, $f(x)$ is a smooth expression since it is the
composition of the smooth atom \verb|multiply| with two smooth expressions.
(The second argument of \verb|multiply| is a smooth expression since it is
itself the composition of the smooth atom \verb|inv_pos| with a smooth
expression.)
\item The function $f(x) = |c^T x / (x^T A x) - b|$ with $A \in \symm^n_{++}$
can be expressed as 
\[
\verb|abs(multiply(c @ x, inv_pos(quad_form(x, A))) - b)|.
\]
When expressed this way, $f(x)$ is an L-convex expression since it is the
composition of the NS-convex atom \verb|abs| with a smooth expression. 
\item The function $f(x) = (\|x - a \|_2 - b)^2$ can be expressed as
\[
\verb|square(norm2(x - a) - b)|.
\]
When expressed this way, $f(x)$ is \emph{not} DNLP-compliant since the atom
\verb|square| is not monotone and its argument is not smooth. However, when
we rewrite it as $f(x) = (\sqrt{ \| x - a \|_2^2 } - b)^2$ and express it as 
\[
\verb|square(sqrt(sum_squares(x - a)) - b)|.
\]
then the expression is smooth since it is the composition of the smooth atom
\verb|square| with a smooth expression. (The first term of the argument of
\verb|square| is a smooth expression since it is itself the composition of
the smooth atom \verb|sqrt| with a smooth expression.)
\item The function $f(x) = (\sin x)^2$ can be expressed as
\verb|square(sin(x))|. When expressed this way, $f(x)$ is a
smooth expression since it is the composition of the smooth atom
\verb|square| with a smooth expression.
\item The function $f(x) = |x|^2$ can be expressed as
\verb|square(abs(x))|. When expressed this way, $f(x)$ is
\emph{not} a smooth expression since the atom \verb|abs| is not smooth.
Although $f(x)$ simplifies to the differentiable function $f(x) = x^2$, the
expression as written is not smooth.
\end{itemize}

\paragraph{DNLP objectives and constraints.}
DNLP supports many types of nonconvex objectives and constraints.
\begin{itemize}
\item An avoidance constraint of the form $\|x - a \|_2 \geq r$, where $a \in
\reals^n$ and $r \in \reals_{+}$ are given, can be expressed as
\verb|sum_squares(x - a) >= r ** 2|. This is DNLP-compliant since the left-hand
side is an L-concave expression (as it is a smooth expression) and the
right-hand side is an L-convex expression (as it is constant and thus a smooth
expression).
\item A discretized dynamics constraint of the form $x_{1} = x_0 + s
\cos(\theta)$, where $x_{1}, x_0, s$, and $\theta$ are variables, can be
expressed as \verb|x1 == x0 + multiply(s, cos(theta))|. This is DNLP-compliant
since both sides are smooth expressions.
\item Minimizing an objective function of the form $\|(Ax)^2 - b \|_1$, where
the square is taken elementwise, is DNLP-compliant when the objective is
expressed as the L-convex expression \verb|norm1(square(A @ x) - b)|.
\end{itemize}
Later we will see applications where constraints and objectives of these forms
arise.

\paragraph{A non-DNLP example.} The DNLP ruleset cannot express every function.
As a simple example, consider the function $f(x) = \min \{\|x - a\|_\infty, \|x -
b\|_\infty)\}$ expressed as
\[
\verb|min(norm_inf(x - a), norm_inf(x - b))|.
\]
This expression is not DNLP-compliant and we are not aware of any algebraic
reformulation that makes it DNLP-compliant. 

\subsection{Connection to DCP}
As discussed in previous sections, DNLP is closely related to DCP. We now make
this relationship explicit. 
\begin{quote}
\emph{A DNLP-compliant problem is one that is DCP when all its smooth atoms are
linearized, regardless of the point of linearization.}
\end{quote}
In particular, L-convex, L-concave, and smooth expressions are those that become
convex, concave, and affine, respectively, after linearizing all smooth atoms
they contain, regardless the point of linearization. This justifies the
terminology linearizable-convex.

\section{Canonicalization}
\label{sec:canonicalization}
In this section we describe how problems conforming to DNLP are canonicalized to
a standard NLP form.  Our canonicalization differs from the approach adopted by
most NLP modeling languages, in which the user-specified problem is not
transformed, and automatic differentiation is used to provide derivative oracles
for the objective and constraint functions. In contrast, in DCP-based modeling
systems for convex optimization, the core idea is to perform extensive
transformations of the original problem formulation into a standard \emph{conic
form} \cite{aps2025mosek, Nesterov1994, boyd2004}, which obviates the need
for derivative oracles based on automatic differentiation. This approach
also gracefully handles functions that are nondifferentiable or defined only on
a restricted domain. We adopt a similar approach to canonicalize problems
conforming to DNLP. 

\subsection{The canonical form}
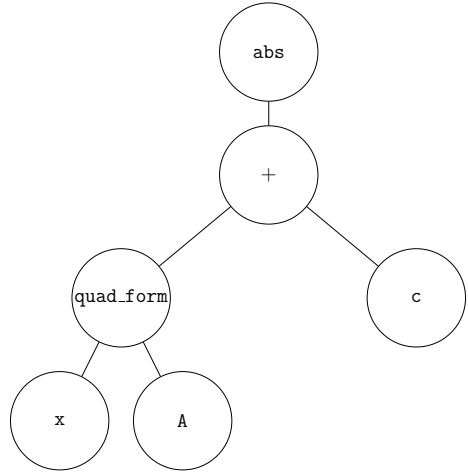
\begin{figure}
\centering
\scalebox{0.65}{ 
\begin{tikzpicture}[
    level distance=25mm,
    level 1/.style={sibling distance=35mm},
    level 2/.style={sibling distance=60mm},
    level 3/.style={sibling distance=25mm},
    every node/.style={circle, draw, minimum size=20mm, inner sep=0pt}
]

\node {\texttt{abs}}
  child {
    node {+}
      child {
          node {\texttt{quad\_form}}
          child { node {\texttt{x}} }
          child { node {\texttt{A}} }
      }
      child {
          node {\texttt{c}}
      }
  };
\end{tikzpicture}
} 
% must have texttt here because of latex compilation
\caption{Expression tree for the L-convex expression $\texttt{abs(quad\_form(x, A) + c)}$.} 
\label{fig:expression-tree}
\end{figure}

The first step of canonicalization is a \emph{parser} that processes the
user-specified problem and constructs one \emph{expression tree} for the
objective and two for each constraint, the left-hand and right-hand sides. 
In an expression tree, each inner node
represents an atom, with its children corresponding to the arguments of the atom.
This is illustrated in figure \ref{fig:expression-tree} for the 
function $f(x) = |x^T A x + c|$ where $A \in \symm_n$ and $c \in \reals$ are
parameters (constants), and $x \in \reals^n$ is the variable, represented by the 
DNLP-compliant expression \verb|abs(quad_form(x, A) + c)|.

Once the expression trees are
constructed, the parser traverses them from the leaves to the root to determine
the smoothness classification of each expression using the definitions given in
\S \ref{sec:dnlp-expressions}. Finally, it verifies that the objective and
constraints conform to the DNLP ruleset described in \S \ref{sec:dnlp-ruleset}. 

After the parser has verified that the problem conforms to DNLP, the
\emph{canonicalizer} traverses the expression trees from the root to the leaves
and transforms the user-specified problem, distinguishing between how 
smooth and nonsmooth atoms are treated.

\paragraph{Smooth atoms.} When a node corresponding to a smooth atom is
encountered, we first check whether the atom has full domain. If not, we
introduce auxiliary variables for its arguments and add constraints linking
these new variables to the original arguments.  We also specify bounds on the
new variables to explicitly encode the domain of the atom. If a smooth atom has
full domain, we apply no transformations to it.%, with two exceptions 
%described in \S \ref{sec:sparsity-encouraging-reformulations}.

A simple example illustrating how smooth atoms are canonicalized is the problem
\[
\begin{array}{ll}
\mbox{minimize} & - \sum_{i=1}^m \log(b_i - a_i^T x)   \\
\mbox{subject to} & \| C x - d \|_2^2 \leq 1,
\end{array}
\]
with variable $x \in \reals^n$. The corresponding canonicalized problem is 
\[
\begin{array}{ll}
\mbox{minimize} & - \sum_{i=1}^m \log(t_i)   \\
\mbox{subject to} & \| C x - d \|_2^2 \leq 1  \\
& t = b - Ax \\
& t \geq 0,
\end{array}
\]
where the variables are (the original one) $x$ and (the new one) $t$. (The
matrix $A \in \reals^{m \times n}$ has rows $a_i^T, \: i =1, \ldots, m$.) Here,
$t$ was introduced for the argument of the logarithm since the log-atom has
restricted domain. No new variable was introduced for the argument to the
squared Euclidean norm in the constraints, since the atom has full domain. Also
note that we explicitly added the bound $t \geq 0$. (Explicitly communicating
function domains via bounds to the solver makes them more robust.) Bounds are
passed to solvers explicitly as bounds, and not as general linear inequality
constraints. This ensures that, as long as the solver respects the variable
bounds strictly (which interior-point solvers do), the solver will never attempt
to evaluate the value or derivatives of an atom at points outside its domain.  

For a problem that only involves smooth atoms, this procedure for traversing the
expression trees results in an equivalent problem formulation similar to a canonical
form proposed by Smith \cite{smith1996Thesis, Smith1996}, known as the
\emph{Smith form}, with the modification that we only introduce
new variables for the arguments of atoms lacking full domain. (In the original
definition of the Smith form, a variable is introduced for any atom argument that is 
not a variable by itself, even if the atom has full domain \cite[table
1]{Smith1996}.) Another distinction from our approach is that the original Smith
form always converts problems into \emph{graph form}, \ie, each nonlinear atom
$\phi$ is replaced by an auxiliary variable $t$ together with the equality
constraint $t = \phi(x)$. For example, instead of minimizing $\phi(x)$ directly,
one minimizes $t$ subject to the constraint $t = \phi(x)$ over $x$ and $t$. 

\paragraph{Nonsmooth atoms.}
When a node corresponding to a nonsmooth atom $\phi$ is encountered, we replace
the atom with an auxiliary variable $t$ and add the constraint $t = \phi(x)$.
Next, we relax this constraint to $t \geq \phi(x)$ if $\phi$ is NS-convex, or to
$t \leq \phi(x)$ if $\phi$ is NS-concave. When the original problem is
DNLP-compliant, this relaxation is (without any further assumptions) guaranteed
to be \emph{lossless} in the sense that (1) the optimal value of the relaxed
problem is the same as that of the original problem, and (2) the set of optimal
$x$-values of the relaxed problem is the same as that of the original problem.
Finally, we express the relaxed constraint using a smooth reformulation, as we
will describe in \S \ref{sec:smooth-epigraph-formulations}.

\subsection{Smooth epigraph formulations}
\label{sec:smooth-epigraph-formulations}
As described in the previous section, any problem conforming to DNLP is
equivalent to a problem in which any atom that is not smooth appears in a
constraint of the form $t \geq \phi(x)$ if $\phi$ is NS-convex, or $t \leq
\phi(x)$ if $\phi$ is NS-concave. The sets $\{(x, t) \mid t \geq \phi(x)\}$ and
$\{(x, t) \mid t \leq \phi(x)\}$ are known as the \emph{epigraph} and
\emph{hypograph} of $\phi$, respectively. Table
\ref{tab:smooth-epigraph-formulations} describes how these are transformed into
smooth formulations that satisfy LICQ. Most of these transformations are
standard and covered in introductory linear programming classes. Automating them
is nevertheless valuable, as the procedure can be tedious and error-prone,
especially when the original problem involves compositions of atoms.

\begin{table}
\centering
\small 
\caption{Smooth epigraph and hypograph formulations of nonsmooth atoms.}
\label{tab:smooth-epigraph-formulations}
\renewcommand{\arraystretch}{1.2} % increase row height
\setlength{\tabcolsep}{5pt} % adjust column spacing
\begin{tabular}{>{\raggedright\arraybackslash}p{2.6cm}
>{\raggedright\arraybackslash}p{4.1cm}
>{\raggedright\arraybackslash}p{2.3cm}
>{\raggedright\arraybackslash}p{5.3cm}}
\toprule
\textbf{Atom} & \textbf{Definition} & \textbf{Smoothness} & \textbf{Epigraph / Hypograph Implementation} \\
\midrule
\verb|abs| & $\phi(x) = |x|$ & NS-convex & Epigraph: $-t \leq x \leq t$ \\ 
\verb|max| & $\phi(x) = \max\{x, y\}$ & NS-convex & Epigraph: $x \leq t, \: y \leq t$ \\
\verb|norm1| & $\phi(x) = \|x \|_1$ & NS-convex & Epigraph: $-v \leq x \leq v, \: \ones^T v \leq t$ \\
\verb|norm2| & $\phi(x) = \|x \|_2$ & NS-convex & Epigraph: $\verb|quad_over_lin|(x, t) - t \leq 0$ \\
\verb|norm_inf| & $\phi(x) = \|x \|_\infty$ & NS-convex & Epigraph: $-t \ones \leq x \leq t \ones$  \\
\verb|huber| & see table \ref{tab:smoothness-classification} & NS-convex & Epigraph: $w^2 + 2M|v| \leq t, \: w + v = x$ \\
\verb|sum_largest| & $\phi(x; k) = \sum_{i=1}^k x_{[i]}$ & NS-convex & Epigraph: \cite[Exercise 5.19]{boyd2004} \\
\verb|min| & $\phi(x) = \min\{x, y\}$ & NS-concave & Hypograph: $x \geq t, \: y \geq t$ \\
\verb|sum_smallest| & $\phi(x; k) = \sum_{i=n-k+1}^n x_{[i]}$ & NS-concave & Hypograph: \cite[Exercise 5.19]{boyd2004} \\
\bottomrule
\end{tabular}
\end{table}

For every atom in table \ref{tab:smooth-epigraph-formulations}, the smooth
reformulation is equivalent to the original epigraph or hypograph constraint,
with one exception. Specifically, for the \verb|norm2| atom, the point $(x, t)
= (0, 0)$ belongs to the epigraph but does not satisfy the smooth reformulation,
since the domain of \verb|quad_over_lin| is $t > 0$ (see table
\ref{tab:smoothness-classification}). Thus, the smooth reformulation excludes
this single point from the feasible set.

Conceptually, this exclusion closely parallels the behavior of interior-point
methods for conic convex optimization such as MOSEK \cite{mosekSolver},
which represent the epigraph of the \texttt{norm2} atom via a second-order cone
constraint. These solvers use barrier functions that enforce strict feasibility
with respect to the cone, ensuring that iterates remain in the cone interior and
thus never reach the origin.

%\subsection{Sparsity-encouraging problem transformations}
%\label{sec:sparsity-encouraging-reformulations}
%
%XXX log-sum-exp, geo-mean, maybe also mention what we have not implemented yet

\subsection{Two advantages of our canonicalization}
\label{sec:advantages-of-canonicalization}

\paragraph{Initialization.} 
One advantage of our canonicalization procedure is that it simplifies the
task of specifying an initial point when atoms have restricted domains. Our
canonical form ensures that the argument to any atom with a restricted domain is
a variable $t$ that appears only as an argument to that atom and in a constraint
of the form $t = f(x)$, where $f(x)$ is an arbitrary expression. Since solvers
require an initial point that lies within the domain of all objective and
constraint functions, we can simply initialize each $t$ within the domain of its
corresponding atom, without ensuring that the constraint $t = f(x)$ holds
initially. This is straightforward to implement by providing an atom-specific
oracle that returns a default initial value within the atom’s domain. All of
this is automated and handled internally, so the user does not need to worry
about it. (If a good starting point for the original variables is known, the
user should of course specify it manually. In this case, we propagate the
starting point to the auxilliary variables introduced during canonicalization by
evaluating the expressions defining them at the given starting point for the
original variables.)

Without this approach, the task of finding an initial point in the intersection
of the domains of the objective and constraint functions falls to the user --- a
task that can be highly nontrivial. For example, consider computing the analytic
center of a polyhedron of the form 
\[
\{x \in \reals^n \: | \: a_i^T x \leq b_i, ~
i = 1, \ldots, m\},
\]
which can be done by minimizing the function 
\[
f(x) = - \sum_{i=1}^n \log(b_i - a_i^T x).
\]
This function is convex and smooth, so we expect the problem to be readily
solved by a solver like Ipopt. For a problem instance where the polyhedron does
not contain the origin, Ipopt crashes in its first iteration when we interface
it using popular NLP modeling languages such as AMPL, GAMS, JuMP, Pyomo, and
CasADi. The reason is that these modeling languages choose the default initial
point to be the origin, which lies outside the domain of the objective function.
In contrast, when we interface Ipopt using our modeling language, it
successfully solves the same problem instance in 14 iterations.

\paragraph{Nonsmoothness.}
Another advantage of our canonicalization procedure is that it seems more
robust for problems involving nonsmooth functions than other NLP modeling
languages that are not based on DNLP. For example, consider the sparse linear
regression problem
\[
\mbox{minimize } \|Ax - b\|_2^2 + \lambda \|x\|_1,
\]
with variable $x$, where $\lambda > 0$ is a regularization parameter. For this
problem we expect the solution to occur at a point of nondifferentiability. When
we specify the problem in our DNLP-based canonicalization pipeline (without
manually reformulating it), Ipopt gracefully solves a random problem instance in
12 iterations. In contrast, when we pass the same formulation to Ipopt via AMPL,
GAMS, JuMP, Pyomo, and CasADi, Ipopt fails to solve the same instance and
terminates after reaching its default maximum number of 3000 iterations. The
issue is that these modeling languages treat the objective as a black box and
supply Ipopt with derivatives via automatic differentiation, even though the
second term is nondifferentiable at the solution.

\subsection{Our implementation}
We have implemented DNLP as an extension to the DCP-based modeling language
CVXPY, available at 
\[
\text{\url{https://github.com/cvxpy/cvxpy}}.
\]
Problems are expressed using standard CVXPY syntax, augmented with smooth
nonconvex and nonconcave atoms including those listed in
table~\ref{tab:smoothness-classification}. (These atoms have previously not been
available in CVXPY, since DCP rules only permit atoms that are either convex or
concave.) For several common atoms we support simpler syntax as a convenience;
for example, squaring all entries of a vector-valued expression \texttt{expr}
can be done using both \texttt{square(expr)} and \texttt{expr ** 2}. 

\paragraph{Some useful functions and features.} The most useful functions and
features of the DNLP extension are summarized below.
\begin{itemize}
\item \verb|problem.is_dnlp()| returns a Boolean indicating whether the problem
is DNLP.
\item \verb|problem.solve(nlp=True)| carries out DNLP canonicalization and
invokes the default NLP solver on the canonicalized problem (assuming the specified
problem is DNLP). The flag \verb|nlp=True| explicitly instructs CVXPY to treat
the problem as a nonlinear program. If omitted, CVXPY attempts to canonicalize
the problem under DCP rules and raises an error if the problem is not
DCP.
\item The \verb|solve()| method accepts the optional keyword argument
\verb|best_of=N|, where \verb|N| is a positive integer. When provided, the
problem is solved \verb|N| times from different random initializations, and the
best solution found is returned. The random starting point for a variable
\verb|x| is drawn uniformly from a user-specified box given by the attribute
\verb|x.sample_bounds|. If \verb|best_of| is used but \verb|x.sample_bounds| is
not provided, no random initialization is done for \verb|x|, unless it
has finite lower and upper bounds and has not been assigned any value.
In that case, the variable is initialized uniformly at random within its bounds.
\item As in CVXPY, the \verb|solve()| method accepts the optional keyword
argument \verb|solver='solver_name'| to specify that the NLP solver
\verb|solver_name| should be used. Directives and options can be passed to the
solver as additional keyword arguments to the \verb|solve()| method.
\item The variable attribute \verb|x.value| can be used to manually set the
initial value for a variable \verb|x|. 
\item As in CVXPY, parameters can be used to specify a family of problems with a
fixed structure. The value of a parameter \verb|p| must be set using
\verb|p.value| before solving the problem, and can be modified between
different solves without recompiling the problem.
\end{itemize}

\paragraph{Supported solvers.}
We currently support the open-source solvers Ipopt
\cite{wachter2006} and Uno \cite{Vanaret2025}, as well as
the commercial solvers Copt \cite{copt} and Knitro \cite{byrd2006}. Knitro
implements several algorithms for nonlinear programming, including an
interior-point method and an augmented Lagrangian method. These can be selected
by specifying the keyword argument \verb|solver='knitro_ipm'| or
\verb|solver='knitro_alm'| in the \verb|solve()| method, respectively. For
example, to use Knitro's interior-point method, one would write 
\begin{quote} 
\hspace{1.75cm} \verb|problem.solve(nlp=True, solver='knitro_ipm')|. 
\end{quote}

\paragraph{Evaluating derivatives.} Modeling languages for NLP must provide
solvers with oracles for evaluating the Jacobian of the constraints and the
Hessian of the Lagrangian. These matrices are typically sparse, and to evaluate
them efficiently we have implemented a differentiation backend in the C
programming language. In our implementation, each expression maintains
knowledge of its derivatives with respect to its arguments. Prior to
optimization, a symbolic preprocessing phase computes the sparsity patterns of
the Jacobian and the Hessian contributions associated with each expression.
During the subsequent numerical phase, executed at each solver iteration, the
expressions evaluate the corresponding numerical derivative values, which are
then assembled into the full Jacobian and Hessian matrices.

When the \verb|best_of| flag is used or parameters are in the problem,
the symbolic preprocessing phase is executed only once,
and the computed sparsity patterns are cached after the
first solve and reused for subsequent solves.

\section{Planning and geometry}
\label{sec:planning-geometry}
In this and the following sections, we present several simple examples
illustrating our DNLP-based modeling language. Most of these can be implemented
in fewer than 10 lines of code, and they are available at
\url{https://github.com/cvxgrp/DNLP-examples}. The examples were solved using
Ipopt, unless otherwise specified.

The code snippets below avoid for-loop constructs where
possible, using vectorized operations instead by specifying axis arguments to
various atoms. This can have a significant performance benefit, so we encourage
users to do so in their own code.

\subsection{Path planning with obstacles}
\paragraph{Problem.}
We seek the shortest path connecting points $a$ and $b$ in $\reals^d$ that
avoids $m$ circles, centered at $p_j$ with radius $r_j$, $j = 1, \dots, m$
\cite{Latombe1991, schulman2014motion}.  After discretizing the
arc-length-parametrized path into a sequence of points $x_0, \ldots, x_n$, the
problem can be written as
\[
\begin{array}{ll}
\mbox{minimize} & L \\
\mbox{subject to} & x_0 = a, \quad x_n = b \\
& \| x_{i+1} - x_i \|_2^2 \leq (L / n)^2, \quad i = 0, \ldots, n-1 \\
& \| x_i - p_j \|_2^2 \geq r_j^2, \quad i = 1, \ldots, n-1, \quad j = 1, \ldots, m \\
& L \geq 0,
\end{array}
\]
with variables $L \in \reals$ and $x_i \in \reals^d, \: i = 0, \ldots, n$. The
problem data are $a \in \reals^d$, $b \in \reals^d$, and $p_j \in \reals^d$ and
$r_j > 0$ for $j = 1, \ldots, m$.  

\paragraph{DNLP specification.}
The code specifying this problem is given below.
\begin{verbatim}
x, L = Variable((d, n + 1)), Variable(nonneg=True)
constr = [x[:, 0] == a, x[:, n] == b,
          sum((x[:, 1:] - x[:, :-1]) ** 2, axis=0) <= (L / n) ** 2]
for i in range(n):
    constr += [sum((x[:, i] - p) ** 2, axis=1) >= r ** 2]
prob = Problem(Minimize(L), constr)
x.value = ...  # initialize to straight line path
prob.solve(nlp=True)
\end{verbatim}

\paragraph{Alternative DNLP-compliant formulations.} The constraint $\| x_{i+1}
- x_i \|_2^2 \leq (L / n)^2$ can also be expressed as $\| x_{i+1} - x_i \|_2
\leq L / n$, which is DNLP-compliant because the left-hand side is L-convex.
Since the objective is decreasing in $L$, these constraints are tight at
optimality, so we can also replace them by equalities of the form $\| x_{i+1} -
x_i \|_2^2 = (L / n)^2$. A constraint of this form is DNLP-compliant, as its
left-hand side is smooth. (Among these three formulations, the first one
converges in the fewest iterations in our experiments.)

\begin{figure}
\centering
\includegraphics[width=0.6\textwidth]{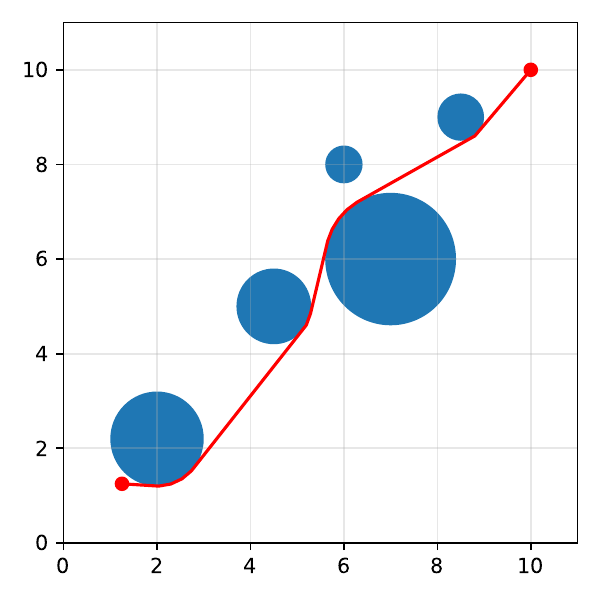}
\caption{Shortest path connecting two points while avoiding circular obstacles.}
\label{fig:path-planning}
\end{figure}

\paragraph{Results.}
We consider a problem instance with dimension $d=2$, $n=50$ path segments, and
$m=5$ obstacles. Figure \ref{fig:path-planning} shows the solution to this
problem instance, when initialized as the straight line path from $a$ to $b$.
For other initializations, the final path is different.

\clearpage 
\subsection{Circle packing}

\paragraph{Problem.}
The goal is to arrange $n$ circles in $\reals^2$ with given radii $r_i$ for $i =
1, \ldots, n$, so that they do not overlap and are contained in the smallest
possible square \cite{Specht2013, Hifi2009}. The optimization problem can be
formulated as 
\[
\begin{array}{ll}
\mbox{minimize} & \max_{i=1,\ldots,n} \left( \| c_i \|_\infty + r_i \right) \\
\mbox{subject to} & \|c_i - c_j \|_2^2 \geq (r_i + r_j)^2, \quad 1 \leq i < j \leq n, \\
\end{array}
\]
where the variables are the centers of the circles $c_i \in \reals^2, \: i = 1,
\ldots, n$, and the radii $r_i$ are given. If $L$ is the value of the objective
function, the circles are contained in the square $[-L, L]^2$.

\paragraph{DNLP specification.}
The code specifying this problem is given below.

\begin{verbatim}
c, constr = Variable((n, 2)), []
for i in range(n - 1):
    constr += [sum((c[i, :] - c[i+1:, :]) ** 2, axis=1) >= 
                   (r[i] + r[i+1:]) ** 2]
cost = max(norm_inf(c, axis=1) + r)
prob = Problem(Minimize(cost), constr)
c.value = uniform(-5.0, 5.0, (n, 2)) # random initial point
prob.solve(nlp=True)
\end{verbatim}

\begin{figure}
\centering
\includegraphics[width=0.6\textwidth]{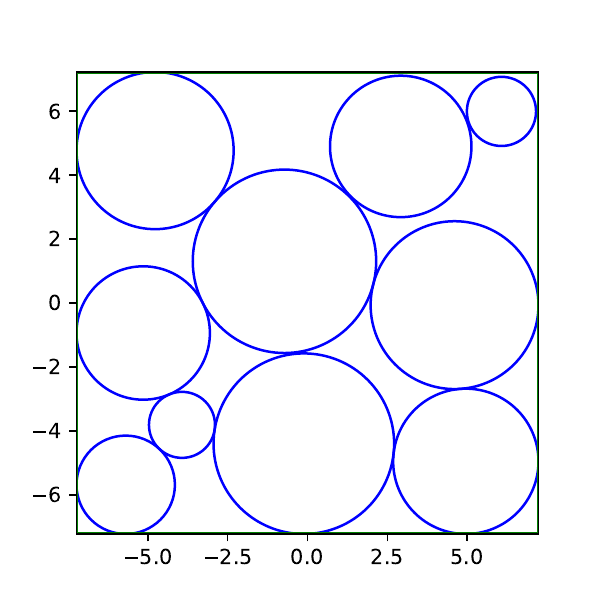}
\caption{Circle packing.}
\label{fig:circle-packing}
\end{figure}

\paragraph{Results.}
We consider a problem instance with $n=10$ circles, with each radius sampled
from a uniform distribution over the interval $[1, 3]$. Figure
\ref{fig:circle-packing} shows one solution to this problem instance, when
initialized with random center locations. The fraction of the square covered by
the circles is 0.72.

\begin{figure}
\centering
\begin{subfigure}{0.45\textwidth}
\centering
\includegraphics[width=\linewidth]{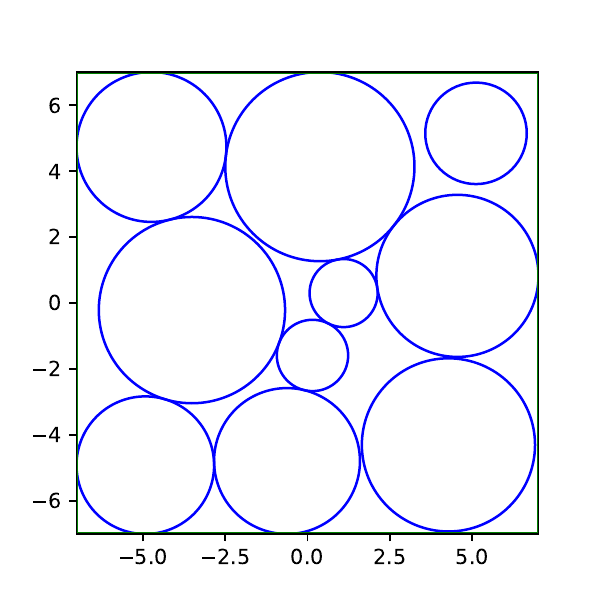}
\end{subfigure}
\hfill
\begin{subfigure}{0.45\textwidth}
\centering
\includegraphics[width=\linewidth]{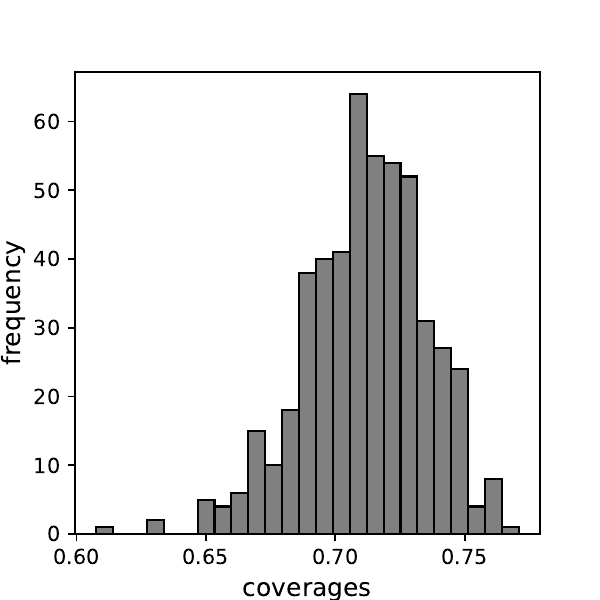}
\end{subfigure}
\caption{The best circle packing found over 500 random initializations (left),
and a histogram of the coverages obtained across all initializations (right).}
\label{fig:circle-packing-best-of}
\end{figure}

%Unlike the previous path planning example where a natural initialization exists
%(the straight line path), this problem has no obvious good initialization. 
To solve the problem multiple times with different random
initializations, we can replace the line \verb|prob.solve(nlp=True)| in the code
snippet above with
\begin{verbatim}
                         
                        c.sample_bounds = [-5.0, 5.0]   
                        prob.solve(nlp=True, best_of=500).
\end{verbatim}
This solves the problem instance 500 times with different random
initializations for the circle centers, each drawn uniformly from the square
$[-5, 5]^2$. With this approach, the fraction of the square covered by the
circles is 0.77 for the best solution found. Figure
\ref{fig:circle-packing-best-of} shows the best solution, along with a histogram
of the coverages obtained across all initializations. 

\clearpage 

\subsection{Nonlinear optimal control}
\paragraph{Problem.}
We consider a simple model of a car in $\reals^2$ as described in \cite[\S
19.4]{Boyd2018}. After time discretization with step size $h > 0$, the state is
$x_k \in \reals^3$, with $((x_k)_1,(x_k)_2)$ denoting its position at time
$t=kh$, and $(x_k)_3$ denoting its angle or orientation. The control input,
which we choose, is $u_k \in \reals^2$, where $(u_k)_1$ is the speed and
$(u_k)_2$ is the steering angle over the time interval $t\in [kh, (k+1)h]$. The
goal is to choose inputs $u_k$ for $k=0,\ldots, N -1$ to move the car from
a given initial state $x^\text{init}$ to a given final state $x^\text{final}$.

The car dynamics are given by $x_{k+1} =  f(x_k,u_k)$, where
\[
f(x_k,u_k) = x_k + (u_k)_1 h
\left[
\begin{array}{c}
\cos(x_k)_3 \\
\sin(x_k)_3 \\
(\tan(u_k)_2) / L
\end{array}
\right]
\]
and $L > 0$ is the wheelbase length of the car. We are given limits $a_{\max}$
and $\omega_{\max}$ on the acceleration and steering angle rate, expressed as $|
(u_{k+1})_1 - (u_k)_1 | \leq a_{\max} h$ and $| (u_{k+1})_2 - (u_k)_2 | \leq
\omega_{\max} h$. We also have lower and upper limits $s_{\min} \leq (u_k)_1
\leq s_{\max}$ and $\phi_{\min} \leq (u_k)_2 \leq \phi_{\max}$ on the speed and
steering angle. We want the control input to be small and smooth, so as
objective we take the sum of the squared Euclidean norms of the control input
over all time steps plus a term that penalizes rapid changes, weighted by
$\gamma > 0$. This gives us the problem
\[
\begin{array}{lll}
\mbox{minimize} & \sum_{k=0}^{N-1} \| u_k \|_2^2 + \gamma \sum_{k=0}^{N-2} \| u_{k+1} - u_k \|_2^2 \\
\mbox{subject to} & x_{k+1} = f(x_k, u_k), & \quad k = 0, \ldots, N-1 \\
& x_0 = x^{\text{init}}, \quad x_N = x^{\text{final}} & \\
& | (u_{k+1})_1 - (u_k)_1 | \leq a_{\max} h, & \quad k = 0, \ldots, N-2 \\
& | (u_{k+1})_2 - (u_k)_2 | \leq \omega_{\max} h, & \quad k = 0, \ldots, N-2 \\
& s_{\min} \leq (u_k)_1 \leq s_{\max}, & \quad k = 0, \ldots, N-1 \\
& \phi_{\min} \leq (u_k)_2 \leq \phi_{\max}, & \quad k = 0, \ldots, N-1,
\end{array}
\]
with variables $x_0, \ldots, x_N$ and $u_0, \ldots, u_{N-1}$. The problem data
are $h$, $L$, $a_{\max}$, $\omega_{\max}$, $s_{\min}$, $s_{\max}$,
$\phi_{\min}$, $\phi_{\max}$, $\gamma$, and the initial and final states
$x^{\text{init}}$, $x^{\text{final}}$.

\paragraph{DNLP specification.}
The code specifying this problem is given below.
\begin{verbatim}
x, u = Variable((N+1, 3)), Variable((N, 2))
cost = sum_squares(u) + gamma * sum_squares(u[1:, :] - u[:-1, :])
constr = [x[0, :] == x_init, x[N, :] == x_final]
constr += [x[1:, :] == x[:-1, :] + h * hstack([
              multiply(u[:, 0], cos(x[:-1, 2])),
              multiply(u[:, 0], sin(x[:-1, 2])),
              multiply(u[:, 0], tan(u[:, 1]) / L)])]
constr += [abs(u[1:, 0] - u[:-1, 0]) <= a_max * h,
           abs(u[1:, 1] - u[:-1, 1]) <= omega_max * h]
constr += [s_min <= u[:, 0], u[:, 0] <= s_max,
           phi_min <= u[:, 1], u[:, 1] <= phi_max]
prob = Problem(Minimize(cost), constr)
prob.solve(nlp=True)
\end{verbatim}
\paragraph{Problem instance.}
We consider a problem instance where the car starts at the origin with zero
orientation, meaning that it is facing right. The final state is $(0.5, 0.5,
-\pi/2)$, \ie, the car should end up half a unit above and to the right
of its starting position, facing down. We use the parameters $L = 0.1$, $N =
50$, $h = 0.1$, and $\gamma = 10$. The acceleration and steering rate limits are
given as $a_{\max} = 0.35$ and $\omega_{\max} = \pi / 10$, and the speed and
steering angle limits are $s_{\min} = -0.15$, $s_{\max} = 0.6$, $\phi_{\min} =
-\pi/8$, and $\phi_{\max} = \pi/8$.

\begin{figure}
\centering
\includegraphics[width=0.8\textwidth]{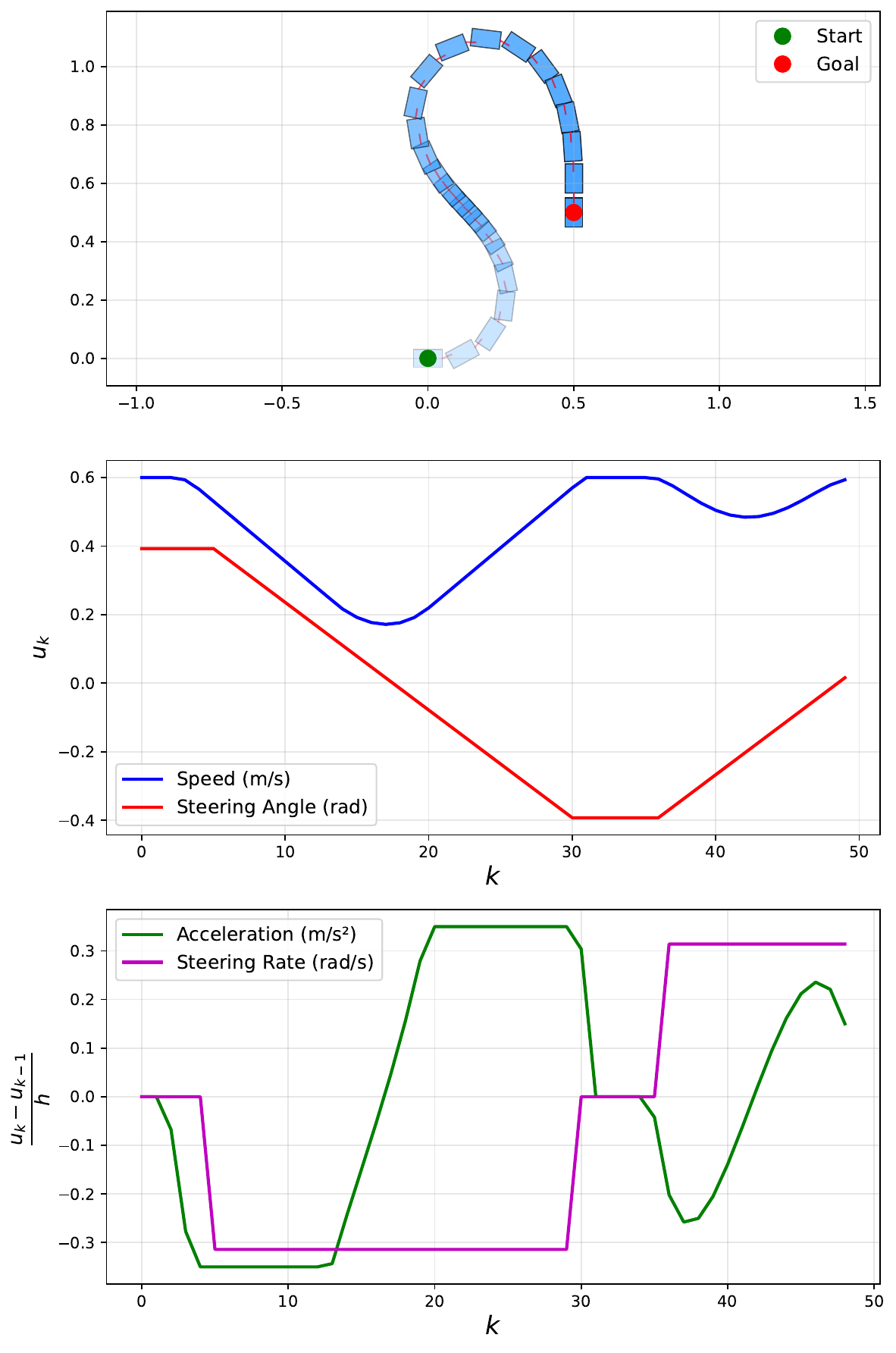}
\caption{Car trajectory. \emph{Top.} Position and orientation of the car.
\emph{Middle.} Speed and steering angle. 
\emph{Bottom.} Acceleration and steering rate.}
\label{fig:car-trajectory}
\end{figure}

\paragraph{Results.}
When we attempted to solve the problem using the default initial point, both
Ipopt and Knitro terminated at an infeasible point. We therefore ran the problem
from 50 random initializations, using the \verb|best_of| flag described in \S
\ref{sec:advantages-of-canonicalization}. Both solvers successfully solved the
problem from some of the initializations. Figure \ref{fig:car-trajectory} shows
the trajectory of the car corresponding to the best solution found, together
with the speed, steering angle, and their rates of change. We see that the
steering angle is initially positive, causing the car to turn left, and then
negative, causing it to turn right, before finally straightening out to reach
the target position facing down.

\clearpage

\section{Signal processing}
\subsection{Location from range measurements}
\paragraph{Problem.}
The goal is to estimate the position of an object from noisy range (distance)
measurements $\rho_i$ to known anchor points $a_i$ in $\reals^2$ for $i = 1,
\ldots, m$ \cite{Smith1987, Beck2008}. We formulate the problem as
\begin{equation} \label{e:location-from-ranges}
\begin{array}{ll}
\mbox{minimize} & \sum_{i=1}^m ( \| x - a_i \|_2 - \rho_i)^2,
\end{array}
\end{equation}
where the variable is the object position $x \in \reals^2$, and the anchor
points $a_i$ and range measurements $\rho_i$ are given. 

\paragraph{DNLP specification.}
The code specifying this problem is given below. To get a DNLP-compliant
formulation, we express $\| x - a_i \|_2$ as $\sqrt{\| x - a_i \|_2^2}$
(see \S \ref{sec:expression-examples}).

\begin{verbatim}
x = Variable(2)
cost = sum_squares(sqrt(sum((x - a) ** 2, axis=1)) - rho)
problem = Problem(Minimize(cost))
problem.solve(nlp=True)
\end{verbatim}

\begin{figure}
\centering
\includegraphics[width=0.95\textwidth]{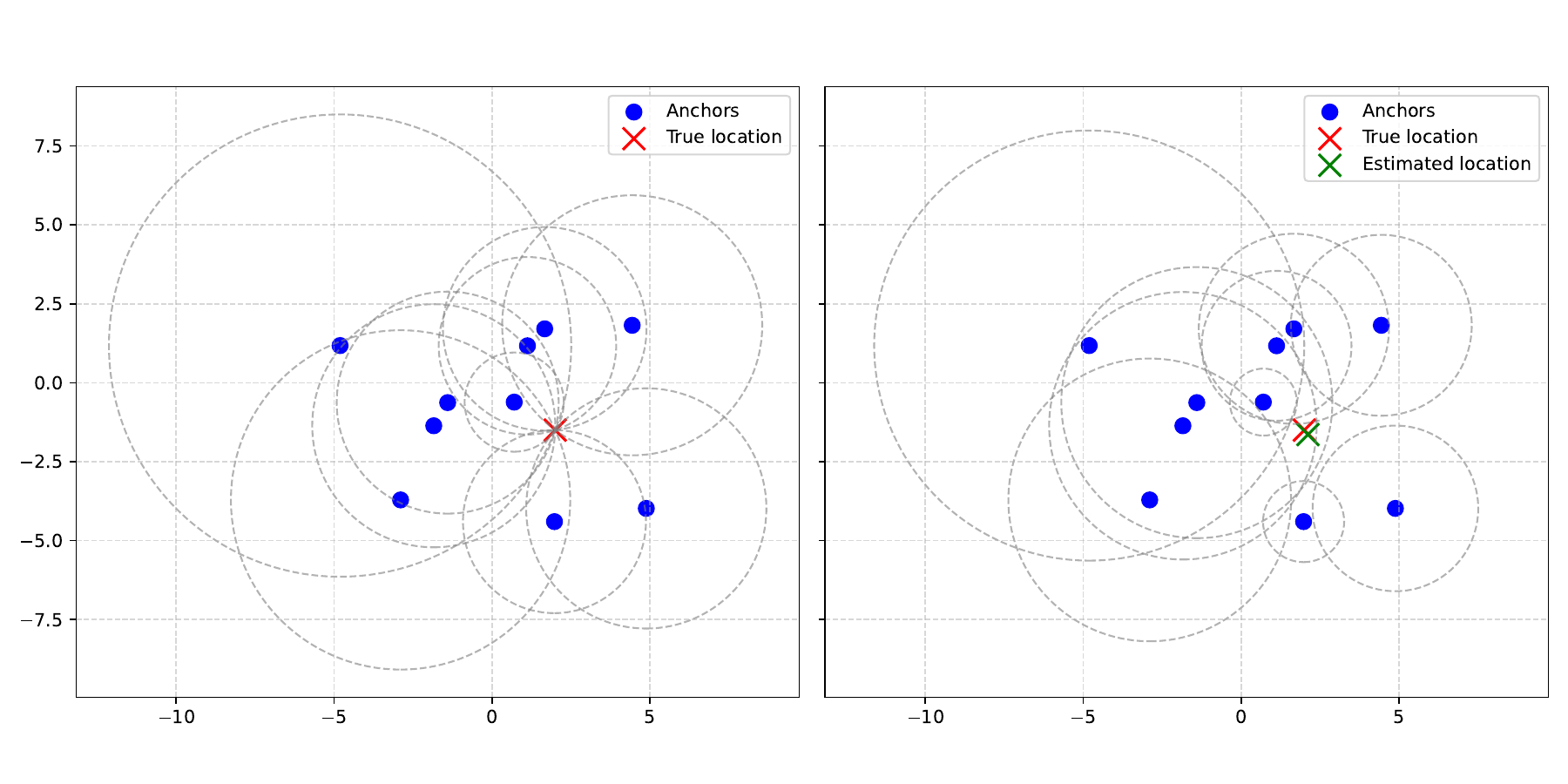}
\caption{Location estimation from range measurements. \emph{Left.}
Range measurements without noise. \emph{Right.} Range measurements with noise.}
\label{fig:location-from-ranges}
\end{figure}

\paragraph{Results.}
We consider a problem instance with $m=10$ anchor points, each sampled from a
uniform distribution over the square $[-5, 5]^2$. We added
zero-mean Gaussian noise with unit standard deviation to the true range
measurements. Figure \ref{fig:location-from-ranges} shows the solution to this
problem instance, with the initial point set to the origin. The dashed
circle around each anchor represents the range measurement from that anchor. The
left figure shows the anchors and range measurements without any noise, and the
true location of the object is at the intersection of the circles. The right
figure shows the noisy range measurements and the estimated location.

\clearpage 

\subsection{Sparse signal recovery}
\paragraph{Problem.}
The goal is to recover a sparse signal $x_0 \in \reals^n$ from a given measurement
vector $y = A x_0$, where $A \in \reals^{m \times n}$ (with $m < n$) is a known
sensing matrix \cite{Candes2008}. A common heuristic based on convex
optimization is to minimize the $\ell_1$ norm of $x$ subject to $Ax = y$. An
alternative approach based on nonconvex optimization is to minimize the sum of
the square roots of the absolute values of the entries of $x$, which tends to
promote sparsity more aggressively \cite{Chartrand2007}. This leads to the problem
\[
\begin{array}{ll}
\mbox{minimize} & \sum_{i=1}^n \sqrt{|x_i|} \\
\mbox{subject to} & Ax = y, \\
\end{array}
\]
with variable $x$. This problem is DNLP-compliant since the objective is L-convex.

\paragraph{DNLP specification.} The code specifying this problem is given below.
For this example, we use Knitro's interior-point method as the solver, because 
Ipopt failed to solve this problem reliably. The issue likely arises from the
fact that the objective function gradient becomes infinite as any entry of $x$
approaches zero, so no KKT point exists for the canonicalized problem.

\begin{verbatim}
x = Variable(n)
cost, constr = sum(sqrt(abs(x))), [A @ x == y]
prob = Problem(Minimize(cost), constr)
prob.solve(nlp=True, solver='knitro_ipm')
\end{verbatim}

\paragraph{Problem instances.} We consider a simulation with signal dimension $n = 100$,
where we vary the number of measurements $m$ from 60 to 80, and the cardinality
of the true signal $x_0$ from 30 to 50. The positions of the nonzero entries of
$x_0$ are sampled from a uniform distribution, with the nonzero values chosen as
$\mathcal{N}(0, 25)$ random variables. The entries of $A$ are sampled from a
standard normal distribution. We say that the recovery is successful if the
relative error $\| \hat{x} - x_0 \|_2 / \| x_0 \|_2$ is less than $10^{-2}$,
where $\hat{x}$ is the recovered signal. To estimate the probability of
successful recovery for each pair of number of measurements and signal
cardinality, we repeat the simulation 100 times and compute the fraction of
successful recoveries. 

\begin{figure}
\centering
\includegraphics[width=\textwidth]{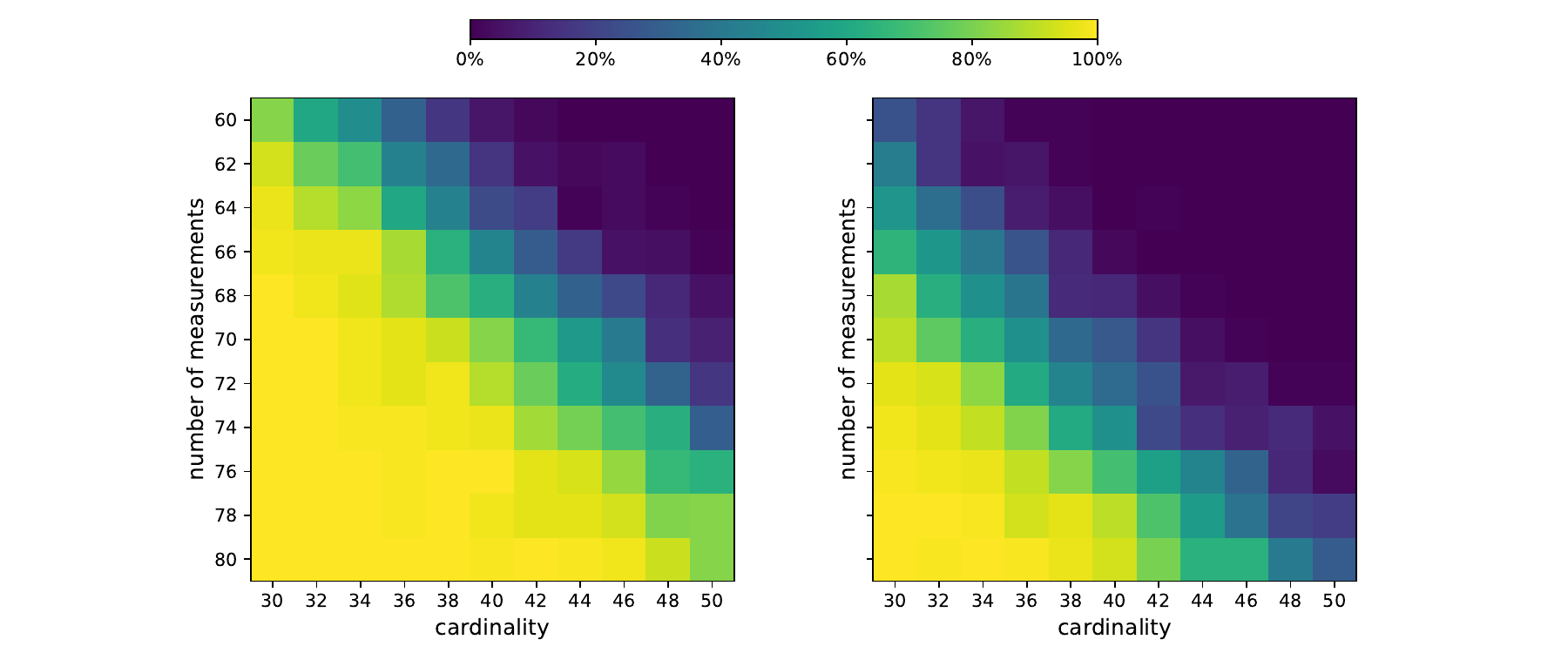}
\caption{Probability of successful signal recovery. \emph{Left.} Approach based
on nonconvex optimization. \emph{Right.} Approach based on convex optimization.}
\label{fig:sparse-signal-recovery}
\end{figure}

\paragraph{Results.}
Figure \ref{fig:sparse-signal-recovery} shows a heatmap
of the estimated probability of successful signal recovery. We see that 
the nonconvex approach is more effective than the convex approach.

\clearpage

\subsection{Phase retrieval}
\paragraph{Problem.}
The goal is to recover a signal $x \in \complex^n$ from the magnitudes of the
complex inner products $a_k^H x$, $k = 1, \ldots, m$, where $a_k \in \complex^n$
are given measurement vectors \cite{Fienup1982, Candes2015}. One version of the
recovery problem can be formulated as 
\[
\begin{array}{ll}
\mbox{minimize} & \||A x|^2 - y^2 \|_1,
\end{array}
\] 
with variable $x \in \complex^n$. Here, $A \in \complex^{m \times n}$ has rows
$a_k^H$, and the absolute value and square operations are applied elementwise.
Since $|Ax|$ is the same if all entries of $x$ are multiplied by a complex
number with unit magnitude, we can only recover $x$ up to some constant phase
shift.

Our current DNLP extension of CVXPY does not support complex variables,
but we can manually reformulate the problem in terms of the real variable
$\tilde{x} = (\Re(x), \Im(x)) \in \reals^{2n}$ as 
\[
\begin{array}{ll}
\mbox{minimize} & \| (B \tilde{x})^2 + (C \tilde{x})^2 - y^2 \|_1,
\end{array}
\]
where the problem data are 
\[
B = \begin{bmatrix}\Re(A) & -\Im(A) \end{bmatrix}
\in \reals^{m \times 2n}, \qquad 
C = \begin{bmatrix} \Im(A) & \Re(A) \end{bmatrix}
\in \reals^{m \times 2n}. 
\]
(Here $\Re(\cdot)$ and $\Im(\cdot)$ denote the real and
imaginary parts, respectively.)

\paragraph{DNLP specification.}
The code specifying this problem is given below.
\begin{verbatim}
x_tilde = Variable(2 * n)
cost = norm1((B @ x_tilde) ** 2 + (C @ x_tilde) ** 2 - y ** 2)
prob = Problem(Minimize(cost))
prob.solve(nlp=True)
\end{verbatim}

\begin{figure}
\centering
\includegraphics[width=0.8\textwidth]{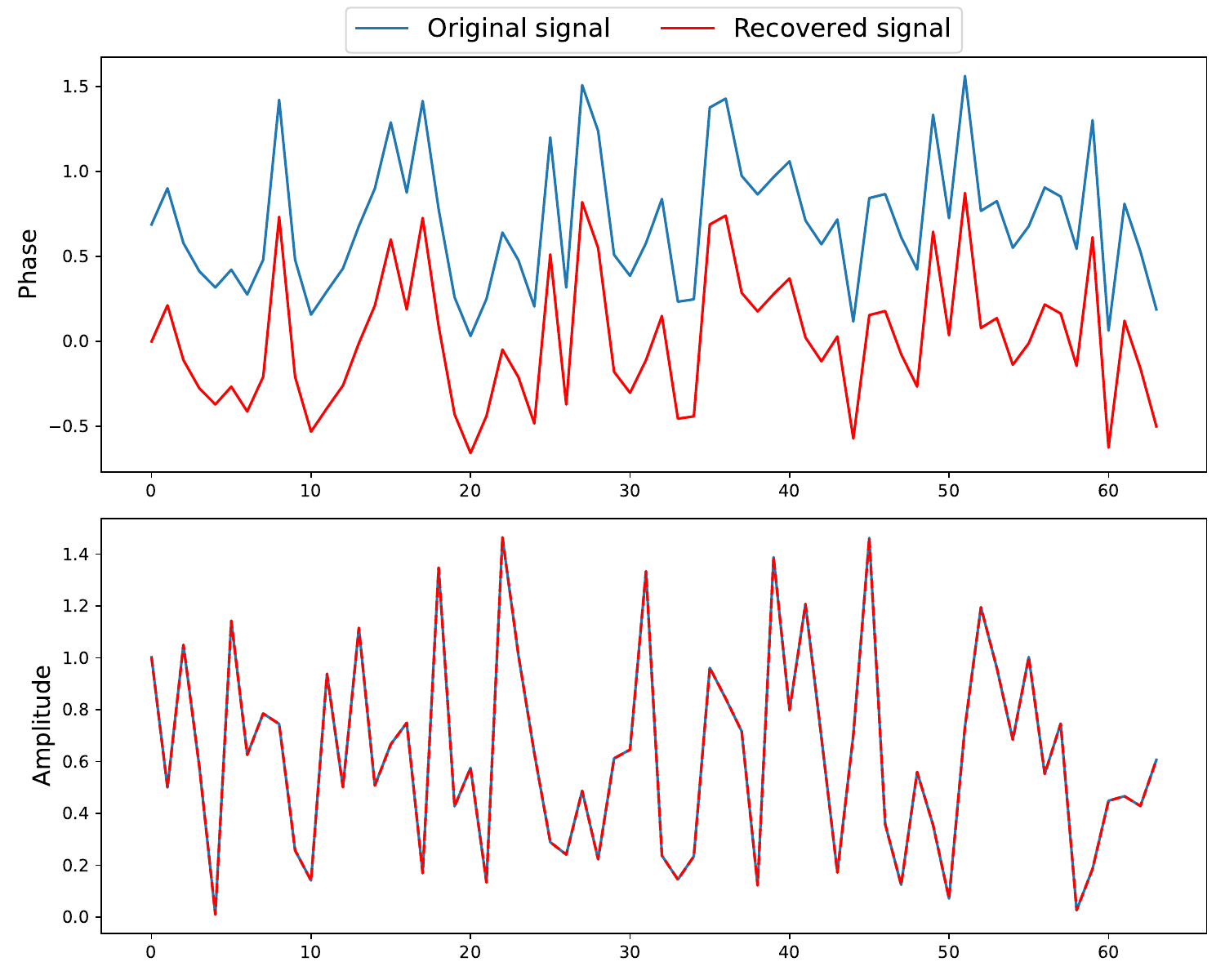}
\caption{Phase retrieval.}
\label{fig:phase-retrieval}
\end{figure}

\paragraph{Results.}
We consider a problem instance with $n = 64$ and $m = 3n$. The real and
imaginary part of each entry of the true signal and the measurement vectors are
sampled uniformly from the unit interval. Figure \ref{fig:phase-retrieval} shows
the original and recovered signals. We see that the signal is accurately
recovered (up to a phase shift).

\clearpage
\section{Finance}

\subsection{Risk-budgeted portfolio construction}
\label{sec:risk-budgeted-portfolio}
\paragraph{Problem.}
Risk-budgeted portfolio construction aims to build a portfolio in which
different sectors contribute specified proportions to the total portfolio risk
\cite{Roncalli2013, Feng2015}. We consider a portfolio of $n$ assets grouped
into $K$ sectors, where $\mathcal{G}_k$ is the set of asset indices in sector
$k$. We let $w_i \geq 0$ denote the fraction of the total portfolio value
(assumed positive) invested in asset $i$. The total portfolio risk is the
standard deviation of the portfolio return $\sigma = (w^T \Sigma w)^{1/2}$,
where $\Sigma \in \symm^n_{++}$ is the asset return covariance matrix. We can
decompose the risk $\sigma$ into components $\sigma_k$ attributable to the
sectors as
\[
\sigma = \frac{w^T \Sigma w}{(w^T \Sigma w)^{1/2}} = 
%\sum_{i=1}^n \frac{w_i (\Sigma w)_i}{(w^T \Sigma w)^{1/2}} = 
\sum_{k=1}^K \sum_{i \in \mathcal{G}_k} \frac{w_i (\Sigma w)_i}{(w^T \Sigma w)^{1/2}} =
\sum_{k=1}^K \sigma_k,
\]
with
\[
\sigma_k = \sum_{i \in \mathcal{G}_k} 
\frac{w_i (\Sigma w)_i}{(w^T \Sigma w)^{1/2}}.
\]
The risk-adjusted return of the portfolio is given by $\mu^T w - \lambda w^T
\Sigma w$, where $\mu$ is the asset return mean, and $\lambda>0$ is a given risk
aversion parameter.

In risk-budgeted portfolio construction, we seek portfolio weights $w$ that
maximize risk-adjusted return subject to sector risks being close to given
proportions $b_k \in (0, 1)$ of the total portfolio risk, \ie, $\sigma_k \approx
b_k \sigma$ for $k = 1, \ldots, K$. With a 10\% tolerance for sector risk
targets, this can be written as the problem
\[\begin{array}{ll} \mbox{maximize} & \mu^T w - \lambda w^T \Sigma w \\
\mbox{subject to} & | \sum_{i \in \mathcal{G}_k} w_i (\Sigma w)_i - b_k  w^T
\Sigma w | \leq 0.1 b_k  w^T \Sigma w, \quad k = 1, \ldots, K \\
& \ones^T w = 1, \quad w \geq 0, \end{array}\] with variable $w \in \reals^n$.

\paragraph{DNLP specification.}
The code specifying this problem is given below. For further efficiency we have
introduced two auxiliary variables $t_1$ and $t_2$ to represent the
subexpressions $\Sigma w$ and $w^T \Sigma w$ that appear multiple times in the
formulation.

\begin{verbatim}
w, t1, t2 = Variable((n, ), nonneg=True), Variable((n, )), Variable()
obj = mu.T @ w - lmbda * t2
constr = [sum(w) == 1, t1 == Sigma @ w, t2 == quad_form(w, Sigma)]
for k, g in enumerate(groups):
    constr += [abs(sum(multiply(w[g], t1[g])) - b[k] * t2)
                                         <= 0.1 * b[k] * t2]
w.value = np.ones(n) / n # uniform initial guess
prob = Problem(Maximize(obj), constr)
prob.solve(nlp=True)
\end{verbatim}

\paragraph{Problem instance.}
We consider a problem instance with $n = 319$ assets from S\&P 500 grouped into
the $K = 5$ largest sectors according to the Global Industry Classification
Standard (GICS), which are Information Technology, Health Care, Financials,
Consumer Discretionary, and Communication Services. The risk budgets are set to
$b = (0.30, 0.25, 0.20, 0.15, 0.10)$, allocating approximately 30\% of portfolio
risk to Information Technology, with the remaining sectors contributing
approximately 25\%, 20\%, 15\%, and 10\%, respectively. We set the covariance
matrix and asset return mean to the sample covariance and empirical mean of
the asset returns, respectively, over the period from January 1, 2020 to January
1, 2025. (Of course, in practice one would use sophisticated methods to estimate
these.) 

\begin{figure}
\centering
\includegraphics[width=0.7\textwidth]{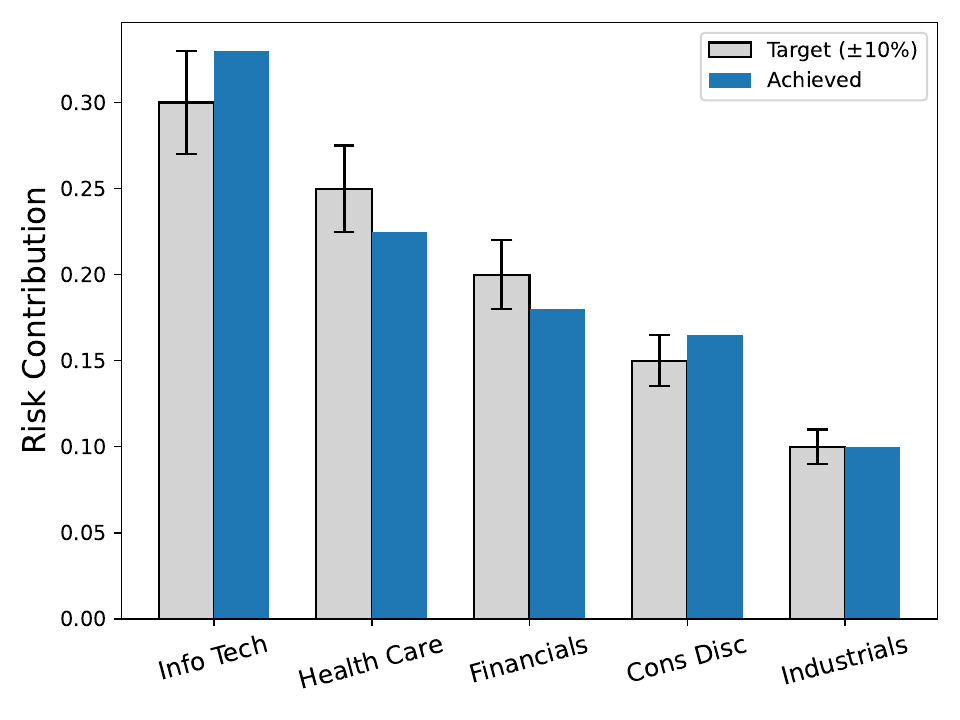}
\caption{Sector risk contributions of the risk-budgeted portfolio.}
\label{fig:risk-budgeted-portfolio}
\end{figure}

\paragraph{Results.} Figure \ref{fig:risk-budgeted-portfolio} shows the sector
risk contributions of the optimized portfolio. 
Two of them take on the highest allowed risk, two take on the smallest
allowed risk, and one is in between the sector risk limits.

\clearpage
\subsection{Optimal product pricing}
\paragraph{Problem.}
We consider the problem of choosing prices of a set of $n$ products to maximize
profit, taking demand responses to price changes into account
\cite{Schaller2025}. Let $p_i^{\text{nom}}$ and $p_i$ denote the nominal and 
new (positive) prices of product $i = 1, \ldots, n$. Likewise, let $d_i^{\text{nom}}$ and
$d_i$ denote the nominal and new (positive) demand of product $i$.
The profit generated by product $i$ is $d_i (p_i - c_i)$, where $c_i$ denotes
the unit production cost. The total profit is thus given as $
P = \sum_{i=1}^n d_i (p_i - c_i)$.

To model the impact of price adjustments on demand, we introduce the logarithmic
relative price change $\pi_i = \log(p_i / p_i^{\text{nom}})$ and the
corresponding logarithmic relative demand change $\delta_i = \log(d_i /
d_i^{\text{nom}})$. We assume a linear relationship between these quantities of
the form $\delta = E \pi$, where $E \in \reals^{n \times n}$ is a given price
elasticity matrix \cite[\S 2F]{mas1995microeconomic}. The total profit can be
expressed in terms of $\delta$ and $\pi$ as
\[
P = \sum_{i=1}^n d_i^{\text{nom}} e^{\delta_i} \left( p_i^{\text{nom}} e^{\pi_i} - c_i \right)
= \sum_{i=1}^n (r_i^{\text{nom}} e^{\delta_i + \pi_i} - \kappa_i^{\text{nom}} e^{\delta_i}),
\]
where $r_i^{\text{nom}} = d_i^{\text{nom}} p_i^{\text{nom}}$ is the nominal
revenue and $\kappa_i^{\text{nom}} = d_i^{\text{nom}} c_i$ is the nominal cost
of product $i$. The problem of choosing prices to maximize profit can thus be
formulated as
\[
\begin{array}{ll} \mbox{maximize} & 
\sum_{i=1}^n (r_i^{\text{nom}} e^{\delta_i + \pi_i} - \kappa_i^{\text{nom}} e^{\delta_i}) \\
\mbox{subject to} & \delta = E \pi, 
\qquad \pi^{\min} \leq \pi \leq \pi^{\max},
\end{array}
\]
with variables $\delta \in \reals^n$ and $\pi \in \reals^n$. The problem data are 
the nominal revenues $r^{\text{nom}} \in \reals^n$, nominal costs $\kappa^{\text{nom}} \in
\reals^n$, price elasticity matrix $E \in \reals^{n \times n}$, and lower and upper bounds
$\pi^{\min}, \pi^{\max} \in \reals^n$ on the relative price changes. 

\paragraph{DNLP specification.} The code specifying the optimal product pricing
problem is given below.
\begin{verbatim}
delta = Variable((n, 1))
pi = Variable((n, 1), bounds=[pi_min, pi_max])
profit = sum(multiply(r_nom, exp(delta + pi)) -
             multiply(kappa_nom, exp(delta)))
constr = [delta == E @ pi]
problem = Problem(Maximize(profit), constr)
problem.solve(nlp=True)
\end{verbatim}

\paragraph{Problem instance.}
We consider a problem instance with $n = 1000$ products. To restrict prices to be
within 10\% of the nominal prices, we set $\pi_i^{\min} = \log(0.9)$ and
$\pi_i^{\max} = \log(1.1)$ for $i = 1, \ldots, n$. We generate the
elasticity matrix $E$ and the nominal revenue and cost vectors $r^{\text{nom}}$ and
$\kappa^{\text{nom}}$ using the same approach as \cite[\S 5]{Schaller2025}.

\paragraph{Results.}
The optimal pricing solution results in a total profit of \$618, corresponding
to a 38\% increase over the \$448 profit achieved with nominal prices.

\clearpage
\subsection{Implied volatility model calibration}
\paragraph{Options.}
European options are financial derivatives that give the holder the right to buy
or sell an underlying asset, such as a stock or index, at a fixed \emph{strike
price} $K$ at a fixed future time $T$, called the \emph{expiry}. Option prices
are quoted both in currency units (\eg, USD) and in terms of their \emph{implied
volatilities} $\sigma^{\text{imp}} > 0$. The implied volatility of an option is
the volatility of the underlying asset for which the Black--Scholes option
pricing formula \cite{Black1973} matches the observed market price of the
option. The \emph{forward price} $F$ is the price agreed today for delivery of
the underlying asset at time $T$. The \emph{log-moneyness} of an option with
strike $K$ is defined as $k = \log(K / F)$. Plotting the implied volatility
against log-moneyness yields the \emph{implied volatility curve}. For more
background on options and implied volatility, see \cite{Hull2021}.

\paragraph{Problem.}
We are given a set of $N$ strikes $K_1, \ldots, K_N$ and the corresponding
implied volatilities $\sigma_1^{\text{imp}}, \ldots, \sigma_N^{\text{imp}}$ for
options on a single underlying asset with fixed expiry $T$ and forward price
$F$. We wish to calibrate a parametric model of the implied volatility curve,
using the observed data. One model for the implied volatility curve is the
\emph{stochastic volatility inspired} (SVI) model \cite{Gatheral2011,
Gatheral2014}, which expresses the implied variance $w =
(\sigma^{\text{imp}})^2$ as a function of the log-moneyness $k$, as
\[
w(k) = \frac{1}{T}\left(a + b \left( \rho (k - m) + 
\sqrt{(k - m)^2 + s} \right)\right),
\]
where $a$, $b$, $\rho$, $m$, and $s$ are model parameters. We seek the
parameters that best fit the observed implied volatilities in the least-squares
sense, subject to known bounds on the model parameters,
\[
\begin{array}{llll}
& a_{\min} \leq a \leq a_{\max}, & b_{\min} \leq b \leq b_{\max}, &
\rho_{\min} \leq \rho \leq \rho_{\max}, \\
& m_{\min} \leq m \leq m_{\max}, & s_{\min} \leq s \leq s_{\max}. & 
\end{array}
\]
This calibration problem can be formulated as
\[\begin{array}{ll} \mbox{minimize} & \sum_{i=1}^N \left(
(\sigma_i^{\text{imp}})^2 - w(k_i) \right)^2 \\
\mbox{subject to} & a_{\min} \leq a \leq a_{\max}, \quad b_{\min} \leq b \leq
b_{\max} \\
& \rho_{\min} \leq \rho \leq \rho_{\max}, \quad m_{\min} \leq m \leq
m_{\max} \\
& s_{\min} \leq s \leq s_{\max}, \end{array}\] with variables $a$, $b$, $\rho$,
$m$, and $s$. The problem data are the implied volatilities
$\sigma_1^{\text{imp}}, \ldots, \sigma_N^{\text{imp}}$, the log-moneyness values
$k_1, \ldots, k_N$, and the variable bounds.

This calibration problem has some convexity properties. It is convex in the
variable $a$, and convex in $b$ when $\rho$, $m$, and $s$ are fixed. Introducing
the new variable $b\rho$, it can be made convex in $a$, $b$, and $b\rho$, when
$m$ and $s$ are fixed. Nevertheless we will formulate and solve it as an NLP.

\paragraph{DNLP specification.} The code specifying this problem is given below.
\begin{verbatim}
a = Variable(bounds=[a_min, a_max])
b = Variable(bounds=[b_min, b_max])
rho = Variable(bounds=[rho_min, rho_max])
m = Variable(bounds=[m_min, m_max])
s = Variable(bounds=[s_min, s_max])

w = (1/T) * (a + b * (rho * (k - m) + sqrt((k - m) ** 2 + s)))
objective = Minimize(sum_squares(w - implied_var))
problem = Problem(objective)
problem.solve(nlp=True)
\end{verbatim}

\paragraph{Problem instance.}
We collect Dow Jones Industral Average options data (DJX) from OptionMetrics
(accessed via Warton Research Data Services) for $N = 30$ options quoted on June
4, 2025, expiring on August 15, 2025. We obtain the risk-free interest rate by
linearly interpolating zero-coupon bonds, and we use put-call parity (see
\cite[\S 18.4]{Hull2021}) to deduce the forward price from the mid call and put prices. We
set the variable bounds to
\[
\begin{array}{llllll}
& a_{\min} = -1, & a_{\max} = 1, & b_{\min} = 0, & b_{\max} = 10 & \rho_{\min} = -1 \\
& \rho_{\max} = 1 & m_{\min} = -2 & m_{\max} = 2, & s_{\min} = 0,  & s_{\max} = 1. 
\end{array}
\]

\begin{figure}
\centering
\includegraphics[width=0.7\textwidth]{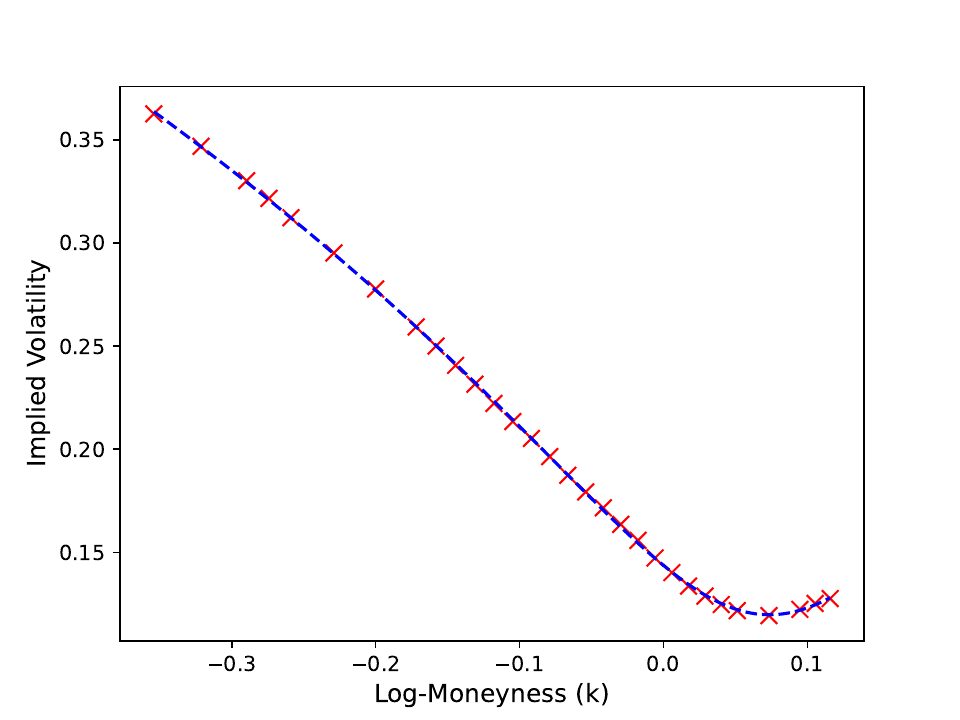}
\caption{SVI model calibration. The red crosses show the observed implied volatilities,
and the blue curve shows the fitted SVI model.}
\label{fig:svi-calibration}
\end{figure}

\paragraph{Results.}
Figure \ref{fig:svi-calibration} shows the observed implied volatilities and the
fitted SVI curve. The fitted parameters are 
\[
a = -0.0083, \quad
b = 0.073, \quad
\rho = -0.068, \quad
m = 0.063, \quad
s = 0.024.
\]

\clearpage
\subsection{Risk-neutral density estimation}

\paragraph{Risk-neutral valuation.}
A key concept in option pricing is \emph{risk-neutral valuation} \cite{Cox1976},
in which the price of a European option is the discounted expected value of its
payoff, where the expectation is taken with respect to a specific
probability distribution called the \emph{risk-neutral distribution}.
Specifically, let $S_T$ denote the price of the underlying asset at maturity $T$
and let $q$ denote the risk-neutral density of $S_T$. 
Then the price (or value) $V$ of a
European option with payoff $f(S_T)$ is
\begin{equation} \label{e:pricing-integral}
V = e^{-rT} \Expect_q[f(S_T)] = e^{-rT} \int_{-\infty}^\infty f(s) q(s) \; ds,
\end{equation}
where $r$ is the risk-free interest rate. In risk-neutral density estimation, we
use observed option prices to fit the density $q$.  This allows us to price
other options for which we do not have a market price, or to detect market
prices that differ from the risk-neutral valuation.

\paragraph{Black--Scholes model.}
In the Black--Scholes model
\cite{Black1973, Black1976}, the terminal asset price is modeled as log-normal,
\[
\log S_T \sim \mathcal{N}(\log F - (1/2)\sigma^2 T, \sigma^2 T),
\]
where $F$ is the forward price of the underlying for delivery at time $T$ and
$\sigma > 0$. The corresponding risk-neutral density is given by
\begin{equation} \label{e:bs-density}
q^{\text{BS}}(s; F, \sigma) = \frac{1}{s \sigma \sqrt{T}}
\phi\left(\frac{\log s - \log F + (1/2)\sigma^2 T}{\sigma\sqrt{T}}\right), 
\end{equation}
where $\phi$ is the probability density function of a Gaussian random variable
with zero mean and unit variance. Under this model, the price of a call option
with strike $K$, which gives the holder the right to buy the underlying at price
$K$ at time $T$ and thus has payoff $f(s) = \max(s - K, 0)$, can be evaluated in
closed form as
\[
C^{\text{BS}}(K; \sigma, F) = e^{-rT} \left(F \Phi(d_1) - K \Phi(d_2)\right),
\]
where $\Phi$ denotes the standard Gaussian cumulative distribution function, and
\[
d_1 = \frac{\log(F/K) + (1/2)\sigma^2 T}{\sigma\sqrt{T}}, \qquad
d_2 = d_1 - \sigma\sqrt{T}.
\]
Similarly, the price of a put option with strike $K$, with
payoff $f(s) = \max(K - s, 0)$, evaluates to
\[
P^{\text{BS}}(K; \sigma, F) = e^{-rT} \left(K \Phi(-d_2) - F \Phi(-d_1)\right).
\]
These celebrated formulas were derived by Black and Scholes
\cite{Black1973, Black1976} and Merton \cite{Merton1973}, work for which
Scholes and Merton were awarded the 1997 Nobel Memorial Prize in Economics,
two years after Black had passed away.

\paragraph{Log-normal mixture model.}
A more flexible model of the risk-neutral density $q$, that can fit observed
prices better, is a mixture of $M$ log-normal components \cite{Richey1990,
Melick1997}. In this model, the risk-neutral density is given by
\[
q(s; F, \sigma, w) = \sum_{j=1}^{M} w_j q^{\text{BS}}(s; F_j, \sigma_j),
\]
where $F \in \reals^M$, $\sigma \in \reals^M_{++}$, and $w \in \reals^M$ are the
parameters of the model, $q^{\text{BS}}$ is the log-normal density
\eqref{e:bs-density}, and $w_j$ are component weights or probabilities, which
must satisfiy $w \geq 0$ and $\ones^T w = 1$. Since the pricing integral
\eqref{e:pricing-integral} is linear in $q$, the call and put prices under this
mixture model are weighted averages of Black--Scholes prices,
\[
\begin{split}
C^{\mathrm{mix}}(K; F, \sigma, w) & = \sum_{j=1}^{M} w_j C^{\mathrm{BS}}(K; \sigma_j, F_j), \\
P^{\mathrm{mix}}(K; F, \sigma, w) & = \sum_{j=1}^{M} w_j P^{\mathrm{BS}}(K; \sigma_j, F_j).
\end{split}
\]

\paragraph{Calibration problem.}
We consider the problem of calibrating the parameters $(F, \sigma, w)$ of the
log-normal mixture model to observed option prices. We are given a set of $N_1$
call strikes $K_1^{\text{call}}, \ldots, K_{N_1}^{\text{call}}$ and the
corresponding observed call prices $C_1, \ldots, C_{N_1}$. We are also given a
set of $N_2$ put strikes $K_1^{\text{put}}, \ldots, K_{N_2}^{\text{put}}$ and
the corresponding observed put prices $P_1, \ldots, P_{N_2}$. We seek the
parameters that best fit the observed call and put prices in the least-squares
sense, subject to the constraint that $q$ is a valid probability distribution
and that the parameters lie within known bounds. This calibration problem can be
formulated as
\[
\begin{array}{ll}
\mbox{minimize} & \sum_{i=1}^{N_1} \omega_i^\text{call}(C_i - C^{\mathrm{mix}}(K_i^{\text{call}}; F, \sigma, w))^2 +
\sum_{i=1}^{N_2} \omega_i^\text{put}(P_i - P^{\mathrm{mix}}(K_i^{\text{put}}; F, \sigma, w))^2 \\
\mbox{subject to} & w \geq 0, \quad \ones^T w = 1 \\
& F_{\min} \leq F_j \leq F_{\max}, \quad \sigma_{\min} \leq \sigma_j \leq \sigma_{\max}, \quad j = 1, \ldots, M, \\
\end{array}
\]
with variables $F \in \reals^M$, $\sigma \in \reals^M$, and $w \in \reals^M$.
The problem data are the observed call and put prices $C_1, \ldots, C_{N_1}$ and
$P_1, \ldots, P_{N_2}$, the strikes $K_1^{\text{call}}, \ldots,
K_{N_1}^{\text{call}}$ and $K_1^{\text{put}}, \ldots, K_{N_2}^{\text{put}}$, the
variable bounds, and the weights $\omega_i^{\text{call}}$ and
$\omega_i^{\text{put}}$ used to give more importance to fitting certain option
prices.

\paragraph{DNLP specification.} The code specifying this problem is given below.

\begin{verbatim}
w = Variable((M, 1), nonneg=True)
sigma = Variable((M, 1), bounds=[sigma_min, sigma_max])
F = Variable((M, 1), bounds=[F_min, F_max])

d1_call = (log(F) -  log(K_call) + 0.5 * sigma ** 2 * T) / (sigma * sqrt(T))
d1_put = (log(F) - log(K_put) + 0.5 * sigma ** 2 * T) / (sigma * sqrt(T))
d2_call = d1_call - sigma * sqrt(T)
d2_put = d1_put - sigma * sqrt(T)

C_BS = multiply(F, normcdf(d1_call)) - multiply(K_call,  normcdf(d2_call))
C_mix = exp(-r * T) * sum(multiply(w, C_BS), axis=0)
P_BS = multiply(K_put, normcdf(-d2_put)) - multiply(F, normcdf(-d1_put))
P_mix = exp(-r * T) * sum(multiply(w, P_BS), axis=0)

obj = Minimize(sum_squares(multiply(sqrt(omega_call), C - C_mix)) +
               sum_squares(multiply(sqrt(omega_put), P - P_mix)))
problem = Problem(objective, [sum(w) == 1])
problem.solve(nlp=True)
\end{verbatim}

\paragraph{Problem instance.}
We use options on the Dow Jones Index (DJX) from OptionMetrics, spanning December
1997 to September 2025. For each standard (non-weekly) expiration date, we select
the cross-section of call and put prices quoted approximately 30 days before  
expiry, discarding expiration dates with fewer than 20 options or missing data.
This yields 334 cross-sections over the sample period, with a median of 56
options per cross-section. We use $M = 3$ mixture components. The volatility
bounds are set to $\sigma_{\min} = 0.05$ and $\sigma_{\max} = 0.7$, and the
forward bounds are set relative to the market-implied forward price
$F^{\mathrm{mkt}}$ as $F_{\min} = 0.7 F^{\mathrm{mkt}}$ and $F_{\max} = 1.3
F^{\mathrm{mkt}}$. (We compute $F^{\mathrm{mkt}}$ using put-call parity.) The
observed option prices $C_i$ and $P_i$ are taken as the mid prices, and the weights
$\omega_i$ are set to the inverse squared bid-ask spread. To assess
out-of-sample pricing accuracy, we use 5-fold cross-validation: the options in
each cross-section are stratified by strike and split into five folds, the model
is fit on four folds, and pricing errors are evaluated on the held-out fold.
Since the problem is nonconvex, we solve each instance from 10 random
initializations (using the \verb|best_of=10|) and keep the best solution.

\begin{figure}
\centering
\includegraphics[width=0.95\textwidth]{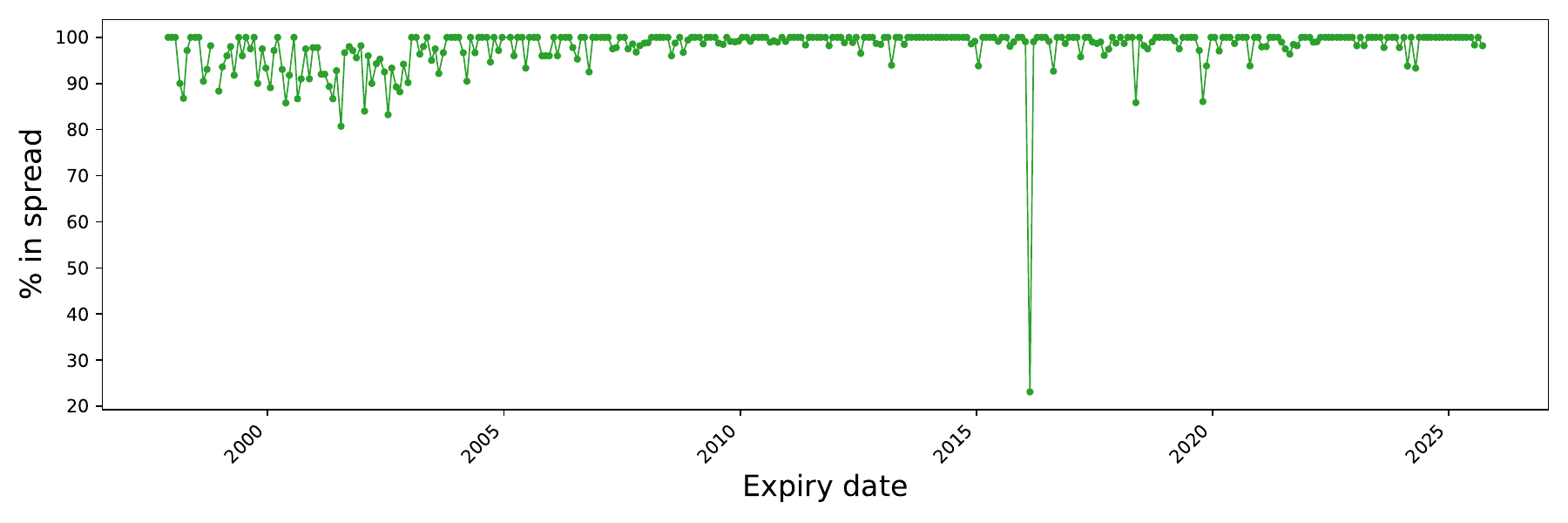}
\caption{The percentage of held-out options for which the model price falls
within the bid-ask spread, across all cross-sections.}
\label{fig:pricing-accuracy}
\end{figure}

\paragraph{Results.}
We measure out-of-sample pricing accuracy using the root-mean-square relative 
pricing error (RMRPE), defined as
\[
\text{RMRPE} = \sqrt{\frac{1}{N}\sum_{i=1}^{N} \left(\frac{\hat{V}_i - V_i}{V_i}\right)^2},
\]
where $\hat{V}_i$ is the predicted mid price and $V_i$ is the observed mid price
of the $i$th held-out option. Over the full sample, the median cross-validated
RMRPE is about 6\%. For comparison, the median relative bid-ask spread across
all options in the dataset is about 8\%, so the typical pricing error is
comparable to the bid-ask spread. As a complementary metric, we report in figure
\ref{fig:pricing-accuracy} the percentage of held-out options for which the
model price falls within the bid-ask spread. 

\clearpage

\section{Energy systems}

\subsection{Optimal power flow}
\label{sec:optimal-power-flow}
\paragraph{Problem.}
The optimal power ﬂow problem seeks to optimize the operation of an electric
power system subject to network power ﬂow constraints and system operating
limits \cite{Frank2016}. We use a standard model, described by a
graph with $n$ buses (nodes), where each bus $i$ is characterized by a voltage
magnitude $v_i$ and a phase angle $\theta_i$. The real and reactive power
injected at the buses are denoted by $p\in \reals^n$ and $q \in \reals^n$,
respectively. These are related to the voltage magnitudes and phase angles via
the equations $p = P \ones$ and $q = Q \ones$, where the bus injection matrices
$P \in \reals^{n \times n}$ and $Q \in \reals^{n \times n}$ are given by 
\begin{equation} \label{e:power-flows}
\begin{split}
P & = (v v^T) \circ ( G \circ C(\theta) + B \circ S(\theta)) \\
Q & = (v v^T) \circ ( G \circ S(\theta) - B \circ C(\theta)).
\end{split}
\end{equation}
Here, $G \in \symm^{n}$ and $B \in \symm^{n}$ are the
(given) real and imaginary parts of the admittance matrix of the network,
$C(\theta) \in \symm^{n}$ and $S(\theta) \in \reals^{n \times n}$ are
defined as
\[\begin{array}{lll}
& C_{ij}(\theta) = \cos(\theta_i - \theta_j), & \quad  S_{ij}(\theta) =
\sin(\theta_i - \theta_j), 
\end{array}\]
and $\circ$ denotes the elementwise (Hadamard) product. Physical limitations of
the network components requires that the power flows and voltages satisfy
certain operational constraints, such as bounds 
\begin{equation} \label{e:operational-constraints}
v^{\min} \leq v \leq v^{\max}, 
\quad p^{\min} \leq p \leq p^{\max}, \quad q^{\min} \leq q \leq q^{\max}.
\end{equation}
(The bounds on $p$ and $q$ can be used to model generation limits at generator
buses and load demands at load buses.) To fix the reference angle of the
network, we force the phase angle at the first bus to be zero. The total
generation cost is typically a convex quadratic function $f(p)$ of the real
power generated at each bus. The optimal power flow problem can thus be
formulated as
\begin{equation} \label{e:prob-optimal-power-flow}
\begin{array}{ll}
\mbox{minimize} & f(p) \\
\mbox{subject to} 
& P  = (v v^T) \circ ( G \circ C(\theta) + B \circ S(\theta)) \\
& Q  = (v v^T) \circ ( G \circ S(\theta) - B \circ C(\theta)) \\
& p = P \ones, \quad q = Q \ones, \quad \theta_1 = 0 \\
& v^{\min} \leq v \leq v^{\max}, \quad p^{\min} \leq p \leq p^{\max}, 
\quad q^{\min} \leq q \leq q^{\max},
\end{array}
\end{equation}
with variables $v$, $\theta$, $P$, $Q$, $p$, and $q$.

\paragraph{DNLP specification.}
The code specifying this problem is given below.
\begin{verbatim}
theta, P, Q = Variable((N, 1)), Variable((N, N)), Variable((N, N))
v = Variable((N, 1), bounds=[v_min, v_max])
p = Variable(N, bounds=[p_min, p_max])
q = Variable(N, bounds=[q_min, q_max])
C, S = cos(theta - theta.T), sin(theta - theta.T)   
constr = [theta[0] == 0, p == sum(P, axis=1), q == sum(Q, axis=1),
          P == multiply(v @ v.T, multiply(G, C) + multiply(B, S)),
          Q == multiply(v @ v.T, multiply(G, S) - multiply(B, C))]
cost = ... # some cost function
prob = Problem(Minimize(cost), constr)
prob.solve(nlp=True)
\end{verbatim}

\paragraph{Alternative DNLP specification.}
The code above declares the bus injection matrices $P$ and $Q$ as dense
matrices, and uses the power flow equations \eqref{e:power-flows} to incorporate
the sparsity pattern of the network only via the admittance matrices $G$ and
$B$. We can also use the variable attribute \verb|sparsity| to explicitly define
$P$ and $Q$ as sparse matrices. If $E$ is the set of edges in the network, we do
this by declaring $P$ and $Q$ as
\begin{verbatim}
                   P = Variable((N, N), sparsity=E)
                   Q = Variable((N, N), sparsity=E).
\end{verbatim}
This alternative approach is more efficient for large networks.

\begin{figure}
\centering
\includegraphics[width=0.7\textwidth]{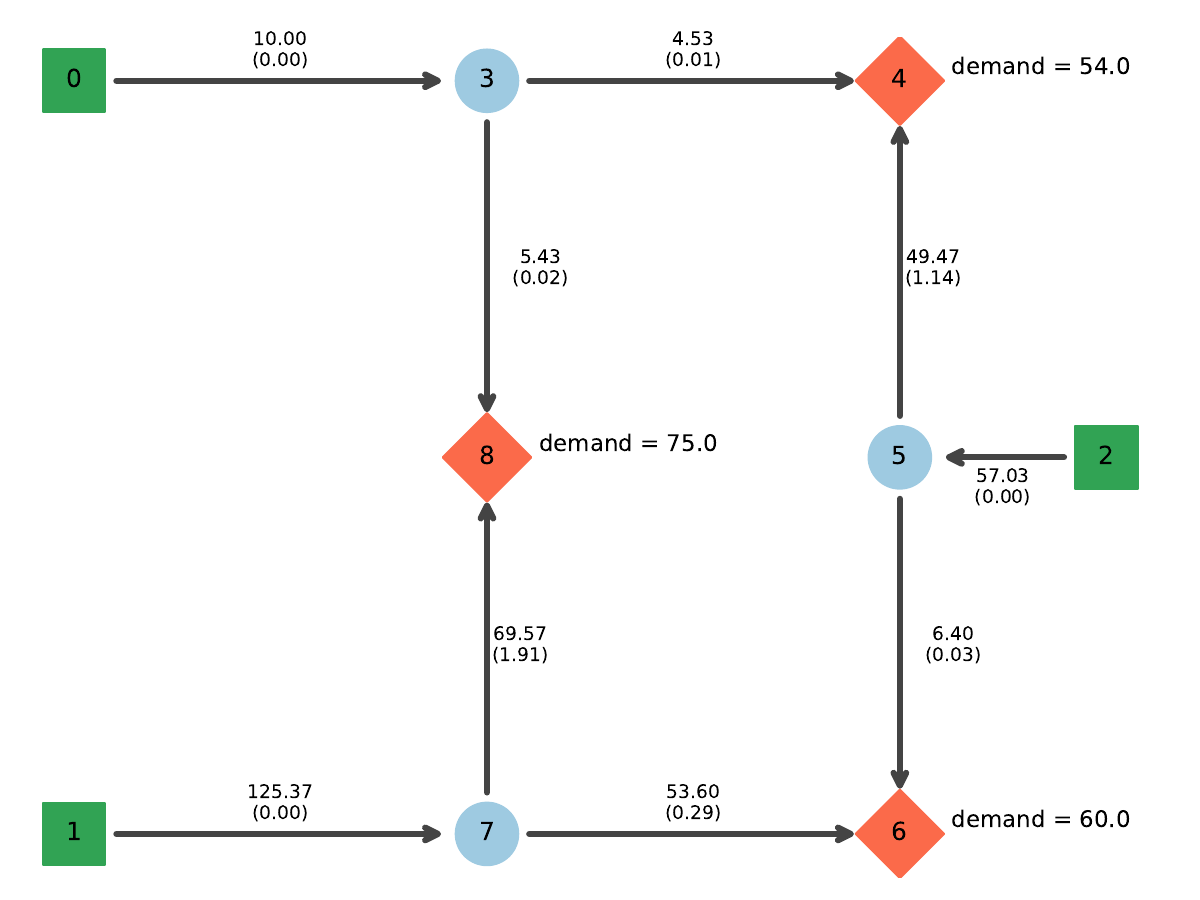}
\caption{Optimal power flow.}
\label{fig:optimal-power-flow}
\end{figure}

\paragraph{Results.}
We consider a 9-node network from \cite{bukhsh2013} with 3
generator buses (green squares), 3 transmission buses (blue circles), and 3 load
buses (orange diamonds). Figure \ref{fig:optimal-power-flow} shows the optimized
real power flow. Each directed edge is annotated with the real power flowing
into the bus at the arrowhead, and, in parentheses, the corresponding real-power
loss on that line. (The real power flow on each line $(i, j)$ is given by
$P_{ij}^{\text{flow}} = P_{ij} - v_i^2 G_{ij}$, with the convention that
positive flow is toward bus $j$. The loss of real power on line $(i, j)$ is
given by $L_{ij} = P_{ij}^{\text{flow}} + P_{ji}^{\text{flow}}$.) The total
generation cost for the computed flow is \$3087.84, which is known to be the
global solution \cite[table 15.2]{Krasko2017}.

\clearpage
\subsection{Power system state estimation}
\paragraph{Problem.}
Closely related to optimal power flow is \emph{power system state estimation}
\cite{Schweppe1970, Abur2004}. Here, the goal is to determine the operating
state (voltage magnitudes and phase angles) of an electric power network from a
set of noisy measurements collected by sensors placed throughout the network. We
consider the same network model as in \S\ref{sec:optimal-power-flow}, with $n$
buses, voltage magnitudes $v \in \reals^n$, phase angles $\theta \in \reals^n$,
and bus power injection matrices $P \in \reals^{n \times n}$ and $Q \in
\reals^{n \times n}$ defined by the power flow equations \eqref{e:power-flows}.
The real and reactive power injections are $p = P\ones$ and $q = Q\ones$.

Sensors placed throughout the network provide $m$ measurements, collected in a
vector $z \in \reals^m$. We model each measurement $z_i$ as a known function $h_i:
\reals^{2n} \to \reals$ of the state $(v, \theta)$, corrupted by additive noise, as
\[
z_i = h_i(v, \theta) + e_i, \quad i = 1, \ldots, m,
\]
where $e_i$ represents random noise with known standard deviation $\sigma_i >
0$. Common measurement types include:
\begin{itemize}
\item Voltage magnitudes at a subset of buses. If measurement $i$
corresponds to a voltage magnitude measurement at bus $k$, then 
$h_i(v, \theta) = v_k$.
\item Real and reactive power injections at buses. For example, if measurement
$i$ corresponds to a real power injection measurement at bus $k$, then $h_i(v,
\theta) = p_k$.
\item Real and reactive line power flows from the sending bus~$i$ to the
receiving bus~$j$ on monitored lines. For example, a real line flow
measurement has $h_i(v, \theta) = P_{ij}^{\text{flow}}$ and a reactive line flow
measurement has $h_i(v, \theta) = Q_{ij}^{\text{flow}}$, where
\[
P_{ij}^{\text{flow}} = P_{ij} - v_i^2 G_{ij}, \qquad
Q_{ij}^{\text{flow}} = Q_{ij} + v_i^2 B_{ij}.
\]
Here, $G_{ij}$ and $B_{ij}$ are the $(i, j)$ entries of the given admittance
matrices $G$ and $B$, and $P_{ij}$ and $Q_{ij}$ are the $(i, j)$ entries of the
bus injection matrices $P$ and $Q$ defined by \eqref{e:power-flows}.
\end{itemize}
Given the measurements $z \in \reals^m$, the standard deviations $\sigma \in
\reals^m_{++}$, and the network admittance data $G \in \reals^{n \times n}$ and
$B \in \reals^{n \times n}$, the estimate of the state is obtained by solving
the weighted least-squares problem
\[
\begin{array}{ll}
\mbox{minimize} & \sum_{i=1}^{m} (1 / \sigma_i^2) (z_i - h_i(v, \theta))^2 \\
\mbox{subject to} & \theta_1 = 0,
\end{array}
\]
with variables $v \in \reals^n$ and $\theta \in \reals^n$. (We have fixed the
phase angle at the first bus to zero to fix the reference angle of the network.)
We treat $P$, $Q$, $p$, and $q$ as functions of $v$ and $\theta$ defined
by the power flow equations \eqref{e:power-flows}, and not as explicit variables
subject to the power flow equations as constraints  (as in the optimal power flow
problem \eqref{e:prob-optimal-power-flow}).
 
\paragraph{DNLP specification.}
The code specifying this problem is given below.
\begin{verbatim}
theta, v = Variable((n, 1)), Variable((n, 1))
C, S = cos(theta - theta.T), sin(theta - theta.T)
P = multiply(v @ v.T, multiply(G, C) + multiply(B, S))
Q = multiply(v @ v.T, multiply(G, S) - multiply(B, C))
p, q = sum(P, axis=1), sum(Q, axis=1)
P_flow = P[fr, to] - multiply(v[fr, 0]**2, G[fr, to])
Q_flow = Q[fr, to] + multiply(v[fr, 0]**2, B[fr, to])
h = hstack([v[V_meas], p[P_meas], q[Q_meas], P_flow, Q_flow])
prob = Problem(Minimize(sum_squares((z - h) / sigma)), [theta[0] == 0])
prob.solve(nlp=True)
\end{verbatim}
Here, \verb|V_meas|, \verb|P_meas|, and \verb|Q_meas| are index arrays
specifying the buses with measurements of voltage magnitudes and real and
reactive power injections, respectively. The arrays \verb|fr| and \verb|to|
specify the sending and receiving buses of the monitored lines.

\paragraph{Results.}
We consider a 3-bus network described in \cite[Example 2.2]{Abur2004}. The
estimated voltage magnitudes are $\hat{v} = (1.000, 0.974, 0.944)$ per unit with
phase angles $\hat{\theta} = (0.0^\circ, {-1.2}^\circ, {-2.7}^\circ)$, which are
the same as the values computed in \cite[Example 2.6]{Abur2004}.

\clearpage

\subsection{Battery model calibration}
\begin{figure}
\centering
\begin{circuitikz}[scale=1, transform shape]
  % OCV source
  \draw (0,0) to[V, v=$v^{\text{oc}}(t)$, invert] (0,3)
  % R0
        to[R, l=$R_0$, -] (3,3)
  % Junction for parallel RC
        -- (4.5,3);
  % R1 (top branch of parallel pair)
  \draw (4.5,3) -- (4.5,3.8)
        to[R, l=$R_1$] (7.5,3.8)
        -- (7.5,3);
  % C1 (bottom branch of parallel pair)
  \draw (4.5,3) -- (4.5,2.2)
        to[C, l=$C_1$] (7.5,2.2)
        -- (7.5,3);
  % U_RC label with +/- circles
  \node[circle, draw, inner sep=1pt, font=\scriptsize] at (4.5,1.2) {$+$};
  \node at (6,1.2) {$U^{\text{RC}}(t)$};
  \node[circle, draw, inner sep=1pt, font=\scriptsize] at (7.5,1.2) {$-$};
  \draw[dashed] (4.5,2.2) -- (4.5,1.5);
  \draw[dashed] (7.5,2.2) -- (7.5,1.5);
  % Output terminal
  \draw (7.5,3) -- (9,3);
  \draw (0,0) -- (9,0);
  % Terminal voltage with +/- circles
  \node[circle, draw, inner sep=1pt, font=\scriptsize] at (9.5,3) {$+$};
  \node[right] at (9.5,1.5) {$v(t)$};
  \node[circle, draw, inner sep=1pt, font=\scriptsize] at (9.5,0) {$-$};
  % Current arrow on main wire
  \draw[->, >=stealth] (8.8,3.3) -- node[above] {$i(t)$} (7.8,3.3);
\end{circuitikz}
\caption{Th\'evenin battery model.}
\label{fig:thevenin-circuit}
\end{figure}
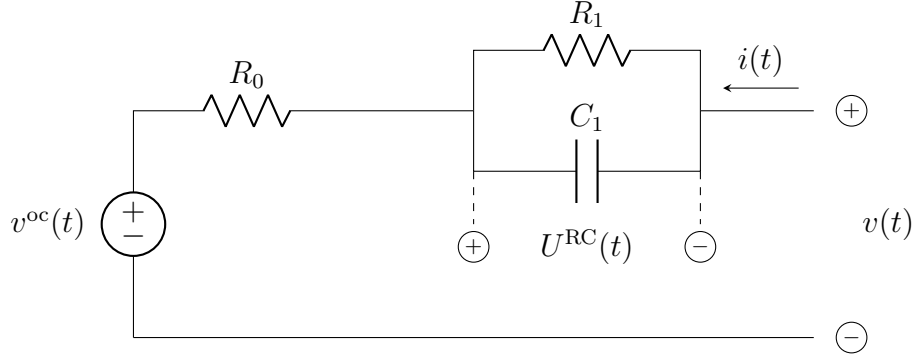

\paragraph{Battery model.}
The goal is to calibrate the unknown parameters of a dynamic model of a
lithium ion (Li-Ion) storage battery, based on experiments.
We consider the so-called Th\'evenin model \cite{He2011} shown in figure
\ref{fig:thevenin-circuit}, which models the battery as an open-circuit voltage
source in series with an internal resistance $R_0$, in $\Omega$ (Ohms)
and a parallel
resistance-capacitance (RC) pair $(R_1, C_1)$, with $R_1$ in $\Omega$ and 
$C_1$ in F (Farads).
We consider a time interval of
$T$ seconds, and denote the charging current at time $t \in [0, T]$ by $i(t)$,
in A (Amperes). We let $q(t)$ denote the stored charge, in C (Coulombs),
$v(t)$ the terminal voltage, in V (Volts), and $v^{\text{oc}}(t)$ the
open-circuit voltage, in V. The terminal voltage is
\[
v(t) = v^{\text{oc}}(t) + R_0 \, i(t) + U^{\text{RC}}(t),
\]
where $U^{\text{RC}}(t)$ is the voltage across the RC pair, which evolves
according to
\begin{equation} \label{e:rc-dynamics}
\frac{dU^{\text{RC}}(t)}{dt} = -\frac{U^{\text{RC}}(t)}{R_1 C_1} + \frac{i(t)}{C_1}.
\end{equation}
The stored charge $q(t)$ satisfies $dq(t)/dt = i(t)$. We model the open-circuit
voltage as a function of the stored charge, as
\[
v^{\text{oc}}(t) = a + \frac{b}{Q^{\text{crit}} - q(t)},
\]
where $a$, $b$, and $Q^{\text{crit}}$ are model parameters. The unit of $a$ is
V, the unit of $b$ is J (Joules), and $Q^{\text{crit}}$ is given in C.
The stored battery charge $q(t)$ is always less than the critical charge 
$Q^\text{crit}$.
This model is parameterized by the six positive parameters
\[
a, \quad b, \quad Q^{\text{crit}}, \quad R_0, \quad R_1, \quad C_1.
\]

\paragraph{Calibration problem.}
We perform an experiment in which we apply a charging current that is constant 
over periods of length $h$ seconds, with $h = T/K$, denoted $i_k$, 
$k=1, \ldots, K$.
We measure the battery terminal voltage at the beginning of each of these
intervals, denoted $v_k^{\text{meas}}$, $k=1, \ldots, K$.
The charge and current are related as
\[
q_{k+1} = q_k + h \, i_k, \quad k = 1, \ldots, K-1,
\]
where $q_k$ is the stored charge at time $t_k$ and we assume that the initial
charge $q_1$ is known. Discretizing the RC dynamics \eqref{e:rc-dynamics} gives
\[
U_{k+1}^{\text{RC}} = \left(1 - \frac{h}{R_1 C_1}\right) U_k^{\text{RC}}
+ \frac{h}{C_1} \, i_k, \quad U_1^{\text{RC}} = 0.
\]
The predicted terminal voltage at time $t_k$ is then $v_k = v_k^{\text{oc}}
+ R_0 \, i_k + U_k^{\text{RC}}$, where $v_k^{\text{oc}} = a + b /
(Q^{\text{crit}} - q_k)$. We wish to find the six parameters $a$, $b$,
$Q^{\text{crit}}$, $R_0$, $R_1$, and $C_1$ that best fit the measured terminal
voltages in the least-squares sense, subject to known bounds:
\[
\begin{array}{llll}
& a_{\min} \leq a \leq a_{\max}, & b_{\min} \leq b \leq b_{\max}, &
Q_{\min}^{\text{crit}} \leq Q^{\text{crit}} \leq Q_{\max}^{\text{crit}} \\
& R_{\min} \leq R_0 \leq R_{\max}, & R_{\min} \leq R_1 \leq R_{\max}, &
C_{\min} \leq C_1 \leq C_{\max}.
\end{array}
\]
This calibration problem can be
formulated as
\[
\begin{array}{lll}
\mbox{minimize} & \sum_{k=1}^{K} (v_k - v_k^{\text{meas}})^2 \\
\mbox{subject to}
& v_k = v_k^{\text{oc}} + R_0 i_k + U_k^{\text{RC}},
& \quad k = 1, \ldots, K \\
& v_k^{\text{oc}} = a + b / (Q^{\text{crit}} - q_k),
& \quad k = 1, \ldots, K \\
& U_{k+1}^{\text{RC}} = \left(1 - \frac{h}{R_1 C_1}\right) U_k^{\text{RC}}
+ \frac{h}{C_1} i_k,
& \quad k = 1, \ldots, K-1 \\
& U_1^{\text{RC}} = 0, \\
& a_{\min} \leq a \leq a_{\max}, \hspace{1.5cm} b_{\min} \leq b \leq b_{\max}, \\
& Q_{\min}^{\text{crit}} \leq Q^{\text{crit}} \leq Q_{\max}^{\text{crit}} 
\hspace{0.8cm} R_{\min} \leq R_0 \leq R_{\max} \\
& R_{\min} \leq R_1 \leq R_{\max}, \hspace{1cm} C_{\min} \leq C_1 \leq C_{\max}.
\end{array}
\]
with variables $v \in \reals^K$, $v^{\text{oc}} \in \reals^K$, $U^{\text{RC}}
\in \reals^K$, and scalar variables $a$, $b$, $Q^{\text{crit}}$, $R_0$, $R_1$,
and $C_1$. The problem data are $v^{\text{meas}} \in \reals^K$,
$i \in \reals^K$, $q \in \reals^K$, $h \in \reals$, and the variable bounds.

\paragraph{DNLP specification.} The code specifying this problem is given below.
\begin{verbatim}
v, v_oc, U_RC = Variable(K), Variable(K), Variable(K)
a = Variable(bounds=[a_min, a_max])
b = Variable(bounds=[b_min, b_max])
Q_crit = Variable(bounds=[Q_crit_min, Q_crit_max])
R0 = Variable(bounds=[R_min, R_max])
R1 = Variable(bounds=[R_min, R_max])
C1 = Variable(bounds=[C_min, C_max])

constrs = [v == v_oc + R0 * i + U_RC, v_oc == a + b / (Q_crit - q),
           U_RC[1:] == (1 - h / (R1 * C1)) * U_RC[:-1] + (h / C1) * i[:-1],
           U_RC[0] == 0.0]
obj = Minimize(sum_squares(v - v_meas))
problem = Problem(obj, constrs)
problem.solve(nlp=True)
\end{verbatim}

\paragraph{Problem instance.}
We generate synthetic data using the discretized model above with 
typical parameters for an 18650 Li-Ion battery,
\[
a = 3.40\mbox{V}, \quad b = 500\mbox{J}, \quad
Q^{\text{crit}} = 6925 \mbox{C},\quad R_0 = 0.10 \Omega, \quad
R_1 = 0.03 \Omega, \quad
C_1 = 1000 \mbox{F},
\]
giving an RC time constant $\tau = R_1 C_1 = 30$ s. The
current profile is a repeating 60-second discharge pattern designed so that
current changes occur on the RC time scale. Specifically, we discharge at $-2$ A
for 40 s, rest at 0 A for 10 s, and charge at $+2$ A for 10 s. This
pattern is repeated 40 times, giving $K = 2400$ time steps with sampling period
$h = 1$ second. The net current per cycle is $-60$ C, so the cell discharges
from the initial charge $q_1 = 4500$ C to approximately 2100 C over 40 minutes.
The measured terminal voltages are generated by adding independent Gaussian
noise with standard deviation $\sigma = 4$ mV to the true voltages. 

\begin{figure}
\centering
\includegraphics[width=\textwidth]{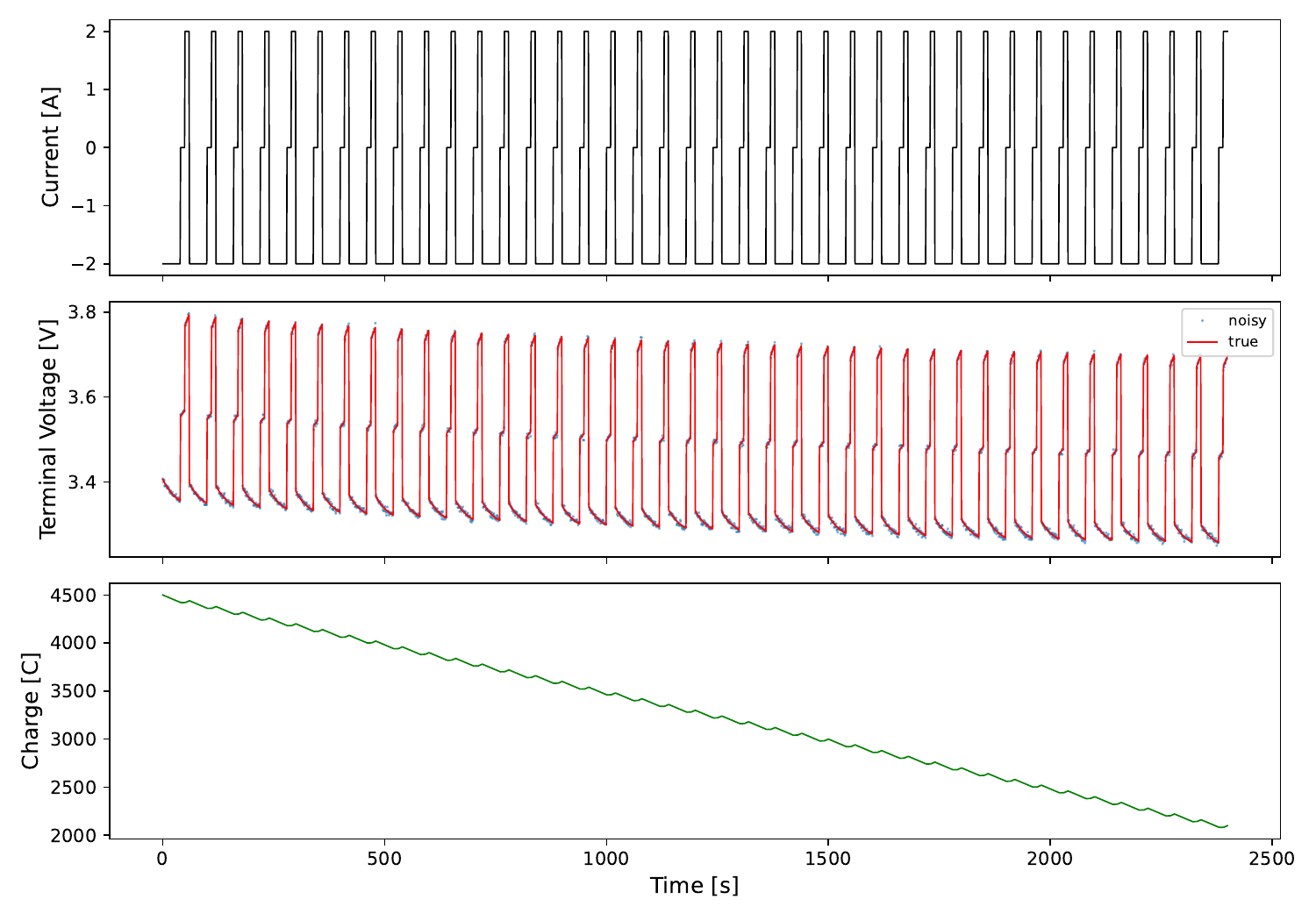}
\caption{Simulated battery data: current profile (top), stored charge (middle),
and measured terminal voltage (bottom).}
\label{fig:battery-data}
\end{figure}

The current
profile, stored charge, and measured terminal voltage are shown in figure
\ref{fig:battery-data}. We set the bounds to 
\[
\begin{array}{lllll}
a_{\min} = 1, & a_{\max} = 10, & b_{\min} = 100, & b_{\max} = 1000, & Q_{\min}^{\text{crit}} = 6000 \\
Q_{\max}^{\text{crit}} = 10000, & R_{\min} = 0.01, & R_{\max} = 0.3, & C_{\min} = 500, & C_{\max} = 2000.
\end{array}
\]

\paragraph{Results.}
The calibrated parameters are 
\[
a = 3.40\mbox{V}, \quad
b = 498\mbox{J}, \quad
Q^{\text{crit}} = 6915 \mbox{C}, \quad 
R_0 = 0.10 \Omega, \quad
R_1 = 0.03 \Omega, \quad
C_1 = 999 \mbox{F}.
\]

\section{Statistics}
\label{sec:statistics}

\subsection{Nonnegative matrix factorization}

\paragraph{Problem.}
The goal is to approximate a given nonnegative matrix $A \in \reals^{m \times
n}$ as the product of two nonnegative matrices $X \in \reals^{m \times k}$ and
$Y \in \reals^{k \times n}$, where $k$ is a given positive integer
\cite{Lee2000, Gillis2020}. One formulation of the problem is
\begin{equation} \label{e:nmf}
\begin{array}{ll}
\mbox{minimize} & \| A - XY \|_F^2 \\
\mbox{subject to} & X \geq 0, \quad Y \geq 0,
\end{array}
\end{equation}
where the variables are the matrices $X$ and $Y$, and $\| \cdot \|_F$ denotes
the Frobenius norm. 

\paragraph{DNLP specification.}
The code specifying this problem is given below.
\begin{verbatim}
X = Variable((m, k), bounds=[0, None])
Y = Variable((k, n), bounds=[0, None])
X.value, Y.value = rand(m, k), rand(k, n)  # random initialization
cost = sum_squares(A - X @ Y)
prob = Problem(Minimize(cost))
prob.solve(nlp=True)
\end{verbatim}

\begin{figure}
\centering
\includegraphics[width=0.95\textwidth]{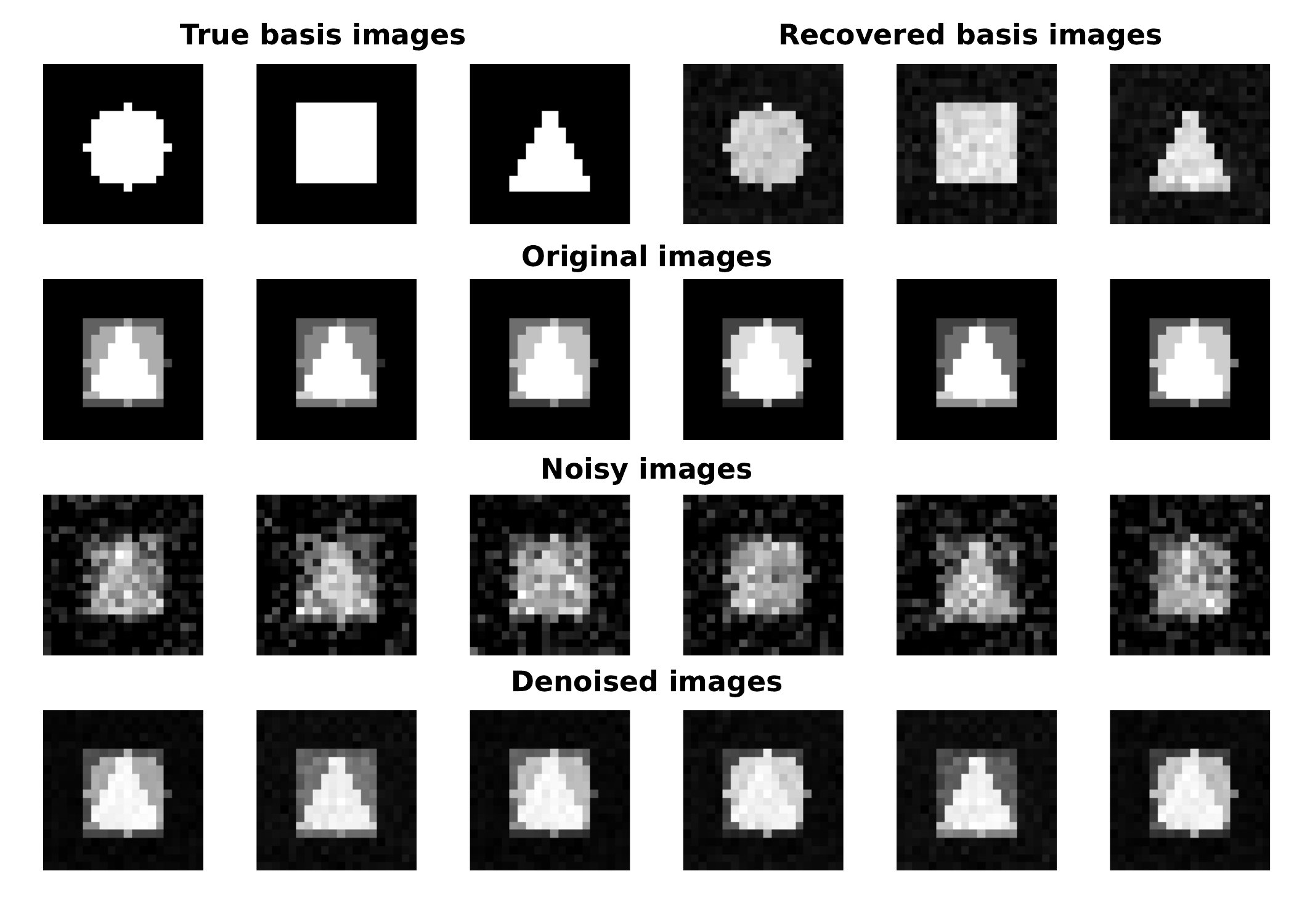}
\caption{Nonnegative matrix factorization for decomposing images into parts.}
\label{fig:nmf}
\end{figure}

\paragraph{Results.}
We use nonnegative matrix factorization to decompose images into basis images
\cite{Lee1999}. First, we generate 100 images of size 20 $\times$ 20 as random
nonnegative combinations of three geometric shapes (a circle, a square, and a
triangle), and then we add noise. After stacking the vectorized noisy images as
columns of a matrix $A \in \reals^{400 \times 100}$, we solve \eqref{e:nmf} with
$k = 3$ to recover the underlying shapes. Figure \ref{fig:nmf} shows the true
basis images followed by the recovered ones (first row), six of the original
images (second row), the same six images after adding noise (third row), and the
denoised images (fourth row) which are given as columns of $X^\star Y^\star$,
where $(X^\star, Y^\star)$ is an approximate solution to \eqref{e:nmf}.

\clearpage 

\subsection{Fitting an exponential decay model}
\paragraph{Problem.}
An \emph{exponential decay model} describes a quantity $y(t)$ that decreases
over time $t$ at a rate proportional to its current value. The model takes the
form
\[
y(t) = a e^{-\lambda t} + c,
\]
where the model parameters are the amplitude $a \in \reals$, the decay rate
$\lambda > 0$, and the asymptotic offset $c \in \reals$. Given noisy
measurements $(t_i, y_i)$ for $i = 1, \ldots, m$, we estimate the parameters by
solving the nonlinear least-squares problem
\[
\begin{array}{ll}
\mbox{minimize} & \sum_{i=1}^m (y_i - a e^{-\lambda t_i} - c)^2 \\
\mbox{subject to} & \lambda \geq 0,
\end{array}
\]
with variables $a, \lambda, c \in \reals$. This problem is not convex, but close
to convex. When $\lambda$ is fixed, it is convex in $a$ and $c$. When $c$ is
known and $a$ is positive, one could also fit $a$ and $\lambda$ by taking
logarithms and solving the linear least-squares problem of minimizing
$\sum_{i=1}^m (\log(y_i - c) - \tilde{a} + \lambda t_i)^2$ subject to $\lambda
\geq 0$, where $\tilde{a} = \log a$ is the log-transformed amplitude.

Alternatively, to gain robustness to outliers, we can replace the least-squares
objective with a Huber loss, resulting in a problem of the form
\[
\begin{array}{ll}
\mbox{minimize} & \sum_{i=1}^m \phi(y_i - a e^{-\lambda t_i} - c; M) \\
\mbox{subject to} & \lambda \geq 0,
\end{array}
\]
where $\phi(\cdot \; ; M)$ is the Huber function with fixed threshold $M > 0$
(see table \ref{tab:smoothness-classification}). Here the variables are $a,
\lambda, c \in \reals$.

Finally, we note that since there are only three variables, it is quite
tractable to solve the problem globally by \emph{gridding}, \ie, evaluating the
objective on a large but finite set of values of the parameters, possibly with
refinement. We do this in our example below to certify the parameters found by
Ipopt as globally optimal. 

\paragraph{DNLP specification.}

The code specifying this problem is given below.
\begin{verbatim}
a, lmbda, c = Variable(), Variable(nonneg=True), Variable()
residuals = y - a * exp(-lmbda * t) - c

# least-squares variant
cost = sum_squares(residuals)
prob = Problem(Minimize(cost))
prob.solve(nlp=True)

# Huber variant
cost = sum(huber(residuals, M))
prob = Problem(Minimize(cost))
prob.solve(nlp=True)
\end{verbatim}

\begin{figure}
\centering
\includegraphics[width=0.7\textwidth]{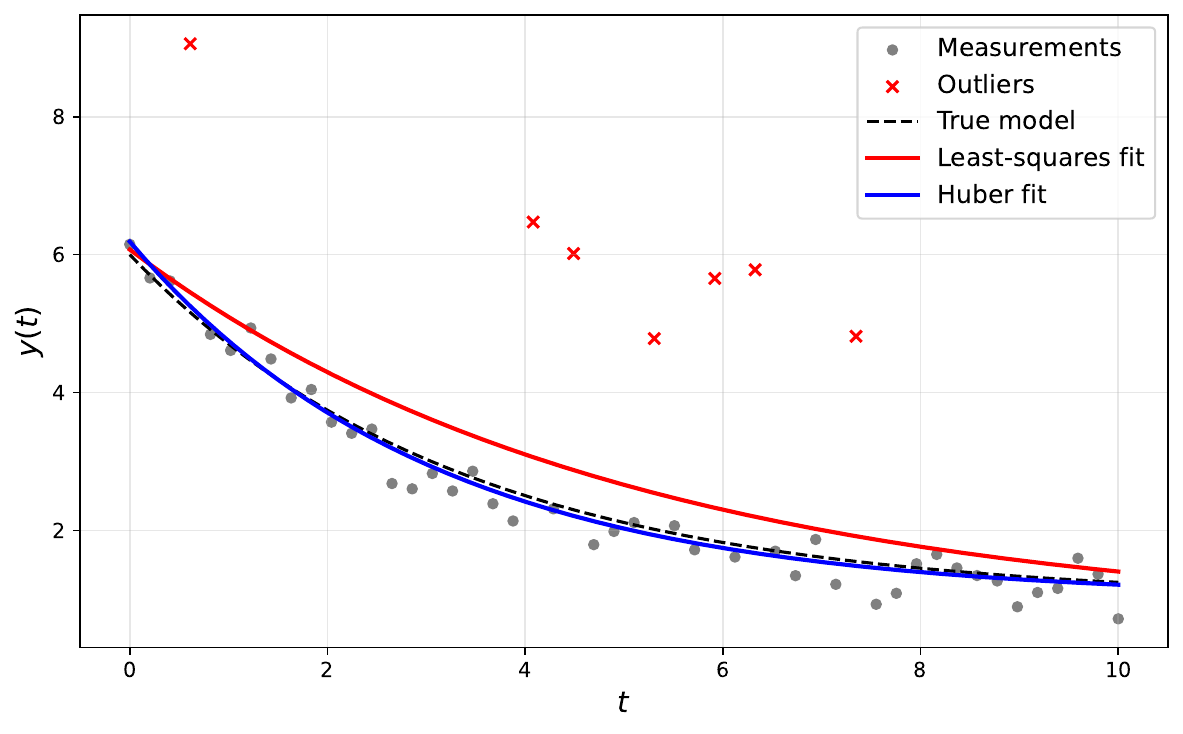}
\caption{Fitting an exponential decay model with outliers.}
\label{fig:exponential_decay}
\end{figure}

\paragraph{Results.}
We generate $m = 50$ noisy measurements from the true model with parameters $a =
5.0$, $\lambda = 0.30$, and $c = 1.0$, adding Gaussian noise. We corrupt $7$ of
the measurements with large outliers. We set the Huber threshold to $M =\sigma$,
where $\sigma$ is the noise standard deviation that is assumed to be known.
Figure \ref{fig:exponential_decay} shows the data, the true model, the
least-squares fit, and the Huber fit. The least-squares fit is visibly pulled
toward the outliers, while the Huber fit remains close to the true model. The
fitted parameters for the least-squares variant and the Huber variant are
$\hat{a} = 5.40$, $\hat{\lambda} = 0.20$, $\hat{c} = 0.67$, and $\hat{a} = 5.2$,
$\hat{\lambda} = 0.33$, $\hat{c} = 1.0$, respectively. The gridding confirms
that Ipopt computes the global minimizers for both variants.

\clearpage

\subsection{Trimmed logistic regression}
\paragraph{Problem.} We are given feature vectors $x_i \in \reals^d$ together
with labels $y_i \in \{-1, 1\}$ for $i = 1, \ldots, N$. We seek a classifier
$\hat{y} = \sign(\theta^T x)$, where $\theta \in \reals^d$ is the model
parameter. In logistic regression, we choose $\theta$ to minimize the logistic
loss
\[
\sum_{i=1}^N \log(1 + \exp(-y_i \theta^T x_i)).
\]
In trimmed logistic regression \cite{Hadi1997, Aravkin2020}, we introduce an
auxiliary weight $w_i \in [0, 1]$ for each data point, allowing the predictor
to downweight outliers and potentially corrupted data points. The parameter
$\theta$ is found by solving
\[
\begin{array}{ll}
\mbox{minimize} & \sum_{i=1}^N w_i \log(1 + \exp(-y_i \theta^T x_i)) \\
\mbox{subject to} & \ones^T w = k, \\
& 0 \leq w_i \leq 1, \quad i = 1, \ldots, N,
\end{array}
\]
with variables $\theta \in \reals^d$ and $w \in \reals^N$. Here, $k \in (0, N)$
is a given parameter that specifies the effective number of samples retained in
the fit.

\paragraph{DNLP specification.} The code specifying the trimmed logistic regression
problem is given below.
\begin{verbatim}
theta = Variable(d)
w = Variable(N, bounds=[0, 1])
loss = sum(multiply(w, logistic(-multiply(y, X @ theta))))
constr = [sum(w) == k]  
prob = Problem(Minimize(loss), constr)
prob.solve(nlp=True)
\end{verbatim}

\paragraph{Problem instance.}
We consider the task of classifying handwritten digits 0 and 1 from the MNIST
data set \cite{Lecun1998}. From the full data set, we randomly select $N = 2000$
images of digits 0 and 1 for training, where each image is represented by $d =
785$ features (the 784 pixel intensities together with an additional bias term).
First, we fit a standard logistic regression model on the clean training data.
We then adversarially corrupt 1\% of the samples by flipping their labels and
refit the standard logistic regression model on this corrupted data. Finally, we
fit a trimmed logistic regression model on the corrupted data using $k = 0.95
N$. To compare the peformance of the different models, we evaluate their
accuracy on a separate test set of 2000 images of digits 0 and 1.

\paragraph{Results.} The standard logistic regression model achieves a test
accuracy of 99.1\% when fitted on the clean training data and 89.3\% when fitted
on the corrupted data. In contrast, the trimmed logistic regression model
achieves a test accuracy of 98.4\% when fitted on the corrupted data.
The weights assigned to the corrupted training samples are zero, indicating that
the trimmed logistic regression model successfully identified and ignored the
corrupted samples.

\clearpage
\subsection{Neural network}
\paragraph{Problem.}
We are given data $x_i \in \reals^{d}$, $y_i \in \{-1, 1\}$, for $i = 1, \ldots,
N$. We seek a classifier $\hat y = \sign(\psi(x; \theta))$, where $\psi(x;
\theta)$ is the output of an $L$-layer neural network with parameters $\theta$.
Hidden layer $k$ contains $n_k$ neurons and is given by
\[
z^{(k)} = \phi(W^{(k)} z^{(k-1)}+ b^{(k)}, \quad k=1,\ldots, L-1,
\]
with $z^{(0)} = x$ as the input layer, where $W^{(k)} \in \reals^{n_{k}\times
n_{k-1}}$, $b^{(k)} \in \reals^{n_{k}}$, $k=1, \ldots, L$, are the weight
matrices and bias vectors to be determined, and $\phi:\reals \to \reals$ is the
activation function, applied elementwise to vectors. The network output is
$\psi(x; \theta) = W^{(L)}z^{(L-1)}+b^{(L)}$, where $\theta = (W^{(1)},b^{(1)},
\ldots, W^{(L)},b^{(L)})$ contains all network parameters. We use $\phi(u) =
\tanh(u)$ as the activation function, since it is smooth. We choose $\theta$ to
minimize the logistic loss plus quadratic regularization,
\[
\frac{1}{N} \sum_{i=1}^N \log(1 + \exp(-y_i \psi(x_i; \theta))) +
\lambda \sum_{k=1}^{L} \|W^{(k)}\|_F^2,
\]
where $\lambda >0$ is a given regularization hyper-parameter, and $\|\cdot\|_F$
is the Frobenius norm. 

\paragraph{DNLP specification.}
The code specifying this problem for a two-layer network is given below.
\begin{verbatim}
W1, b1 = Variable((n1, n0)), Variable((n1, 1))     
W2, b2 = Variable((1, n1)), Variable()   

Z1 = tanh(W1 @ X_train + b1)
psi = W2 @ Z1 + b2
logistic_loss = sum(logistic(-multiply(y_train, psi))) / N_train
regularization = lmbda * (sum_squares(W1) + sum_squares(W2))
objective = Minimize(logistic_loss + regularization)
prob = Problem(objective)

# intialize weights and solve
W1.value = np.random.randn(n1, n0) 
prob.solve(nlp=True)
\end{verbatim}
Here, the training data is stored column-wise, with $X_\text{train} \in
\reals^{n_0 \times N}$ and $y_\text{train} \in \reals^{1 \times N}$.

\paragraph{Problem instance.}
We generate $N = 200$ perturbed samples from two interleaving half-circles in
$\reals^2$, divided into $100$ training and $100$ test samples. The network has
$L = 2$ layers with $n_0 = 2$ input features and $n_1 = 4$ neurons in the hidden
layer. We use three different values of the regularization parameter: $\lambda =
0.0$, $\lambda = 0.001$, and $\lambda = 0.01$. We initialize the weights $W_1$
of the hidden layer randomly. (Without this initialization, all network
parameters are initialized to zero and Ipopt declares that the origin is a local
solution.)

\begin{figure}
\centering
\begin{subfigure}[b]{0.40\textwidth}
\centering
\includegraphics[width=\textwidth]{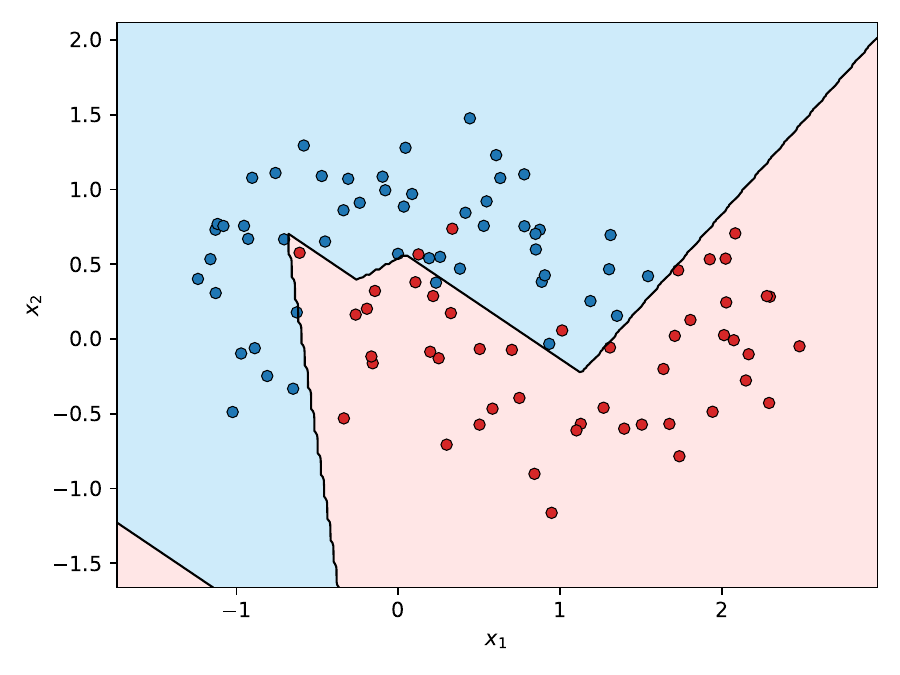}
\caption{$\lambda = 0$}
\end{subfigure}
\hfill
\begin{subfigure}[b]{0.40\textwidth}
\centering
\includegraphics[width=\textwidth]{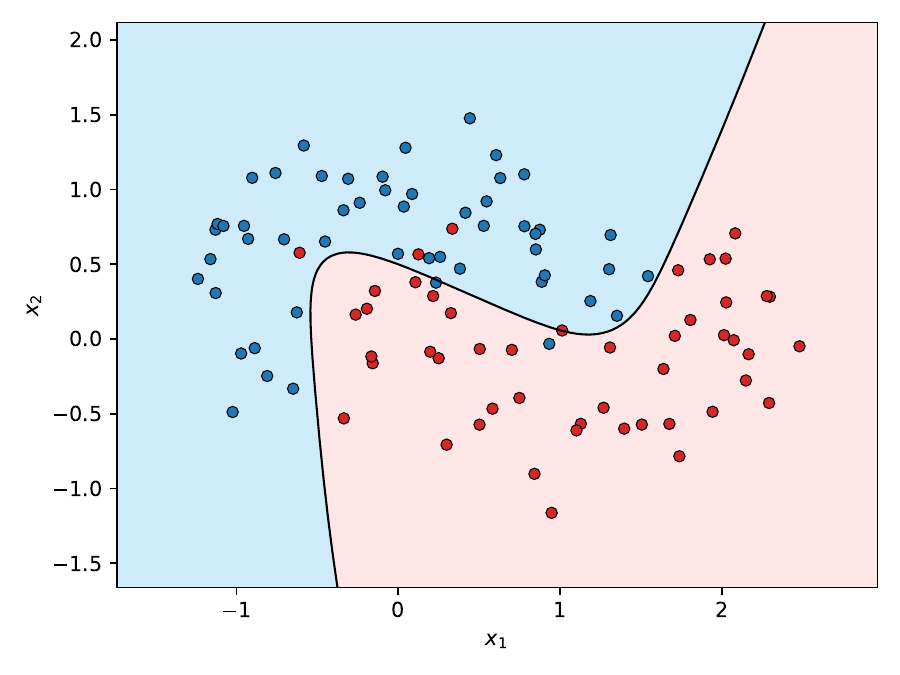}
\caption{$\lambda = 0.001$}
\end{subfigure}
\hfill 
\begin{subfigure}[b]{0.40\textwidth}
\centering
\includegraphics[width=\textwidth]{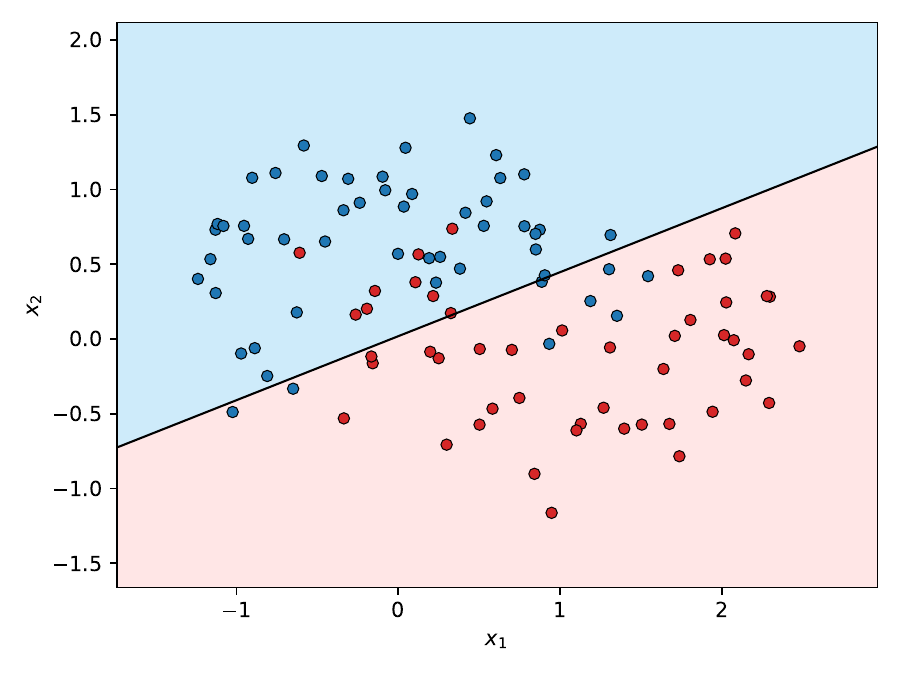}
\caption{$\lambda = 0.01$}
\end{subfigure}
\caption{Training data and decision boundary of neural network classifier for
three values of the regularization parameter $\lambda$.}
\label{fig:nn-classification}
\end{figure}

\paragraph{Results.} Figure \ref{fig:nn-classification} shows the training data
and the decision boundary of the fitted model for the three different values of
$\lambda$. Table \ref{tab:nn-accuracy} shows the training and test accuracies
for the three models. 

\begin{table}[h]
\centering
\caption{Training and test accuracies for neural network classification.}
\label{tab:nn-accuracy}
\begin{tabular}{lcc}
\toprule
 & Training accuracy & Test accuracy \\
\midrule
$\lambda = 0.0$     & 96 \% & 87\% \\
$\lambda = 0.001$   & 94 \% & 94\% \\
$\lambda = 0.01$    & 83 \% & 84\% \\
\bottomrule
\end{tabular}
\end{table}

\clearpage

\subsection{Factor model fitting}
\paragraph{Statistical factor models.}
We are given centered data $x_1, \ldots, x_N \in \reals^d$, assumed to be
independent samples from a zero-mean Gaussian distribution with covariance
matrix $\Sigma \in \symm^d_{++}$. A \emph{statistical factor model} for the
covariance matrix is a factorization of the form
\[
\Sigma = D + F F^T,
\]
where $D \in \reals^{d \times d}$ is diagonal with positive diagonal entries and
$F \in \reals^{d \times r}$, where $r$ is a positive integer typically much
smaller than $d$. Statistical factor models are widely used in finance, signal
processing, and other fields to model high-dimensional data with low intrinsic
dimensionality (see, \eg, \cite{Bartholomew2011, Johansson2023, Cederberg2026}.)
The maximum-likelihood estimation problem for fitting $F$ and $D$ is of the form
\cite{Joreskog1967}
\begin{equation} \label{e:factor-model-ml-problem}
\begin{array}{ll}
\mbox{minimize} & \log \det(D + F F^T) + \Tr((D + F F^T)^{-1} S) 
\end{array}
\end{equation}
with variables $F$ and diagonal $D$, where the problem data is the sample
covariance matrix $S = (1/N) \sum_{i=1}^N x_i x_i^T$. We may also want to
include the constraints $D_{ii} > 0$.

\paragraph{Problem.} The golden standard for solving
\eqref{e:factor-model-ml-problem} is the expectation-maximization (EM) algorithm
\cite{Rubin1982}. However, we can also solve it with a general-purpose NLP
solver. This requires reformulating the problem in a form suitable for such
solvers, which do not natively handle positive definite matrix variables.

First, from the matrix inversion lemma \cite[\S C.4]{boyd2004}, we have
\[
\begin{split}
(FF^T + D)^{-1} & = D^{-1} - D^{-1}F (I + F^T D^{-1} F)^{-1} F^T D^{-1} = E - G G^T, 
\end{split}
\] where $E = D^{-1}$ and $G = D^{-1} F (I + F^T D^{-1} F)^{-1/2}$. This suggests
that \eqref{e:factor-model-ml-problem} is equivalent to 
\begin{equation} \label{e:factor-model-ml-problem-reformulated}
\begin{array}{ll}
\mbox{minimize} & -\log \det(E - G G^T) + \Tr((E - G G^T) S) 
\end{array}
\end{equation}
with variables $G \in \reals^{d \times r}$ and diagonal $E \in \reals^{d \times
d}$. If $(G^\star, E^\star)$ is an optimal solution to this problem, then an
optimal solution to \eqref{e:factor-model-ml-problem} is given by $F^\star =
(E^{\star})^{-1} G^\star (I - (G^{\star})^T (E^{\star})^{-1} G^\star)^{-1/2}$
and $D^\star = (E^{\star})^{-1}$.% \cite{XXX}.

General-purpose NLP solvers can natively handle the second term in the
objective, but not the first term involving the log-determinant. To reformulate
it, we first reduce the dimension of the matrix inside the log-determinant from
$d$ to $r$ via the Schur-complement identity \cite[\S A5.5]{boyd2004}:
\[
\det(E - G G^T) = \det(I - G^T E^{-1} G) \det(E),
\]
where $I \in \reals^{r \times r}$ is the identity matrix. Next, we introduce an
auxiliary lower-triangular variable $L \in \reals^{r \times r}$ together with
the constraint $L L^T = I - G^T E^{-1} G$, which lets us write $\log \det(I -
G^T E^{-1} G) = 2 \sum_{i=1}^r \log L_{ii}$. Finally, it is useful to express
$\Tr((E - G G^T) S)$ as $\diag(S)^T e - \|R^T G \|^2_F$, where $e \in \reals^d$
is the vector of diagonal entries of $E$, $R \in \reals^{d \times d}$ is the
lower triangular Cholesky factor of $S$, and $\| \cdot \|^2_F$ denotes the
Frobenius norm. After these transformations, an equivalent formulation of
\eqref{e:factor-model-ml-problem-reformulated} is
\[
\begin{array}{ll}
\mbox{minimize} & -\sum_{i=1}^d \log(e_i) - 2 \sum_{i=1}^r \log(L_{ii}) 
+ \diag(S)^T e - \|R^T G \|^2_F \\
\mbox{subject to} & L L^T = I - G^T \diag(e)^{-1} G, 
\end{array}
\]
with variables $G$, $L$, and $e$. 

\paragraph{DNLP specification.} The code specifying this problem is given below.

\begin{verbatim}
e = cp.Variable((d, 1), nonneg=True)
L = cp.Variable((r, r), lower_triangular=True)
G = cp.Variable((d, r))
cost = (- sum(log(e)) - 2 * sum(log(diag(L))) + diag(S) @ e 
        - sum_squares(R.T @ G))
constraints = [L @ L.T == eye(r) - G.T @ (G / e)]
prob = cp.Problem(cp.Minimize(cost), constraints)
prob.solve(nlp=True)
\end{verbatim}

\paragraph{Results.}
We consider a problem instance with $d = 100$ and $r = 10$. We plant a
ground-truth factor model by sampling $F_{ij} \sim \mathcal{N}(0, 1)$ and
$D_{ii} \sim \mathcal{U}(0.5, 1.5)$, and then generate $N = 200$ samples from
$\mathcal{N}(0, FF^T + D)$. For Ipopt to converge, we had to specify a
reasonable initial point, which we obtained via principal component factor
analysis \cite[\S 9.3]{Johnson2007}. With this initialization, Ipopt converged
to the same solution as the EM algorithm. 

\clearpage
\section{Conclusions}
\label{sec:conclusions}
In this paper we introduced DNLP, a grammar for specifying nonlinear programs.
Inspired by DCP, DNLP allows smooth functions to be freely combined with
nonsmooth convex and concave functions, provided that a minimal set of rules is
followed. Any problem conforming to DNLP can be canonicalized in a lossless way
to an equivalent smooth NLP, which is then passed to a standard NLP solver. We
described an open-source implementation of DNLP as an extension to CVXPY.

We emphasize that DNLP does not, and cannot, resolve the difficulties that are
intrinsic to a given problem. NLPs are generally nonconvex, and no modeling
language can guarantee that an NLP solver will converge to a global (or even a
local) minimizer. What DNLP \emph{does} address is a more modest but practically
important class of failures: those that arise from na\"{i}ve modeling of
nonsmooth terms, or from initialization issues caused by atoms with restricted
domains. By rewriting such problems into a smooth canonical form and handling
initialization automatically, DNLP removes a common source of avoidable solver
failures. It does not, however, guarantee that every canonicalized problem will
be well behaved.

To demonstrate the breadth of problems that DNLP can express, we surveyed a
range of applications drawn from many different fields. A recurring theme is
that each of these problems can be specified in just a handful of lines of code,
with the canonicalization, initialization, and interface to the solver handled
automatically. We hope that this lowers the barrier to using NLP in practice and
encourages its adoption in new application areas.

\clearpage

\section*{Acknowledgments}
We acknowledge many helpful discussions with several colleagues, including
Maximilian Schaller, Bennet Meyers, Antoine Lesage-Landry, Bartolomeo Stellato,
Clara Baynham, Oscar Dowson, Anthony Degleris, Kevin Tracy, David Pérez Piñeiro,
Philipp Schiele, Art Owen, Michael Salerno, Aleksandr Aravkin, Raphael
Chincilla, Kasper Johansson, and Alexandros Tzikas. We want to thank Amit
Solomon for developing cvxtorch, which was used in an early prototype of this
work. We also want to thank Youssouf Emine (Knitro), Wujian Jack (Copt), Charlie
Vanaret (Uno) for their help with the respective NLP solver interfaces. Finally,
we thank Steven Diamond for very helpful comments and suggestions on the NLP
interface design and its implementation in CVXPY.

William Zhang was supported by Maxime Fortin through the Polytechnique Montreal
Bourse Prestige scholarship.
William Zhang est financ\'e par Maxime Fortin gr\^ace \`a la Bourse
Prestige de Polytechnique Montr\'eal.

\newpage
\bibliographystyle{abbrv}
\bibliography{refs}
\end{document}